\documentclass[11pt, a4paper]{article}
\usepackage{amsmath,amssymb,latexsym,color,enumerate}
\usepackage{authblk}
\usepackage[mathscr]{eucal}
\usepackage[colorlinks,
            linkcolor=blue,
            anchorcolor=green,
            citecolor=magenta
           ]{hyperref}
\newtheorem{thm}{Theorem}[section]

\newtheorem{lemma}[thm]{Lemma}
\newtheorem{cor}[thm]{Corollary}

\newtheorem{conj}{Conjecture}
\newtheorem{obs}[thm]{Observation}

\newtheorem{rema}[thm]{Remark}

\newtheorem{example}{Example}[section]
\newtheorem{defin}{Definition}[section]

\newcommand{\proof}{{\it Proof.\quad}}
\newcommand{\qed}{\hfill\Box\medskip}

\usepackage{CJK}

\begin{document}

\renewcommand{\baselinestretch}{1.3}

\title{\bf $r$-cross $t$-intersecting families for vector spaces}

\author[1]{Mengyu Cao\thanks{E-mail: \texttt{caomengyu@mail.bnu.edu.cn}}}
\author[1]{Mei Lu\thanks{E-mail: \texttt{lumei@tsinghua.edu.cn}}}
\author[2]{Benjian Lv\thanks{Corresponding author. E-mail: \texttt{bjlv@bnu.edu.cn}}}
\author[2]{Kaishun Wang\thanks{E-mail: \texttt{wangks@bnu.edu.cn}}}
\affil[1]{\small Department of Mathematical Sciences, Tsinghua University, Beijing 100084, China}
\affil[2]{\small Laboratory of Mathematics and Complex Systems (Ministry of Education), School of Mathematical Sciences, Beijing Normal University, Beijing 100875, China}
\date{}
\maketitle

\begin{abstract}
Let $V$ be an $n$-dimensional vector space over the finite field $\mathbb{F}_q$, and ${V\brack k}$ denote the family of all $k$-dimensional subspaces of $V$. The families $\mathcal{F}_1\subseteq{V\brack k_1},\mathcal{F}_2\subseteq{V\brack k_2},\ldots,\mathcal{F}_r\subseteq{V\brack k_r}$  are said to be $r$-cross $t$-intersecting if $\dim(F_1\cap F_2\cap\cdots\cap F_r)\geq t$	for all $F_i\in\mathcal{F}_i,\ 1\leq i\leq r.$ The $r$-cross $t$-intersecting  families $\mathcal{F}_1$, $\mathcal{F}_2,\ldots,\mathcal{F}_r$ are said to be non-trivial if $\dim(\cap_{1\leq i\leq r}\cap_{F\in\mathcal{F}_i}F)<t$. In this paper, we first determine the structure of $r$-cross $t$-intersecting families with maximum product of their sizes. As a consequence, we partially prove one of Frankl and Tokushige's conjectures about $r$-cross $1$-intersecting families for vector spaces. Then we describe the structure of non-trivial $r$-cross $t$-intersecting families $\mathcal{F}_1$, $\mathcal{F}_2,\ldots,\mathcal{F}_r$ with maximum product of their sizes under the assumptions
$r=2$ and $\mathcal{F}_1=\mathcal{F}_2=\cdots=\mathcal{F}_r=\mathcal{F}$, respectively, where the $\mathcal{F}$ in the latter assumption is well known as $r$-wise $t$-intersecting family. Meanwhile, stability results for non-trivial $r$-wise $t$-intersecting families are also been proved.
\medskip



\end{abstract}

\section{Introduction}
Let $n$ and $k$ be integers with $1 \le k \le n$, and $V$ an $n$-dimensional vector space over the finite field $\mathbb{F}_q$, where $q$ is necessarily a prime power. We use ${V\brack k}$ to denote the family of all $k$-dimensional subspaces of $V$.  In the sequel we will abbreviate ``$k$-dimensional subspace'' to ``$k$-subspace''. Recall that  for any positive integers $a$ and $b$ the \emph{Gaussian binomial coefficient} is defined by
$$
{a\brack b} = \prod_{0\leq i<b}\frac{q^{a-i}-1}{q^{b-i}-1}.
$$
In addition, we set ${a\brack 0}=1$, and ${a\brack c} =0$ if $c$ is a negative integer. Recall that the size of ${V\brack k}$ is equal to ${n\brack k}$.

For any positive integer $t$, a family $\mathcal{F}\subseteq{V\brack k}$ is called $t$-\emph{intersecting} if $\dim(A\cap B)\geq t$ for all $A,B\in\mathcal{F}$. A $t$-intersecting family $\mathcal{F}\subseteq{V\brack k}$ is called \emph{trivial} if all its members contain a common specified $t$-subspace of $V$ and \emph{non-trivial} otherwise. In general, the triviality of a $t$-intersecting family is determined by the following parameter.
For any $\mathcal{F}\subseteq{V\brack k}$ (not necessary to be $t$-intersecting), a subspace $T$ of $V$ is called to be a \emph{$t$-cover} of $\mathcal{F}$ if $\dim(T\cap F)\geq t$ for all $F\in \mathcal{F}$, and the \emph{$t$-covering number} $\tau_t(\mathcal{F})$ of $\mathcal{F}$ is the minimum dimension of a $t$-cover. It is clear that a $t$-intersecting family $\mathcal{F}$ is trivial if and only if $\tau_t(\mathcal{F})=t$.

The structure of $t$-intersecting families of ${V\brack k}$ with maximum size had been completely determined \cite{Deza-Frankl-1983,PL,PR,Hsieh-1975-1,WN,Tanaka-2006-903}, which are known as the Erd\H{o}s-Ko-Rado Theorem for vector spaces. In \cite{AB}, using the $1$-covering number, the structure of non-trivial $1$-intersecting families of ${V\brack k}$ with maximum size was described. In \cite{Cao-vec}, using the $t$-covering number, the authors described the structure of maximal non-trivial $t$-intersecting families of ${V\brack k}$ with large size, from which the extremal non-trivial $t$-intersecting families are determined (see also \cite{D'haeseleer} for the latter result). Recently, in \cite{Ge-2021}, the authors considered an inverse problem for $t$-intersecting families of subspaces, and  provided structural characterizations for the families with maximal total intersection number.

Let $t,$ $r,$ $k_1,\ldots, k_r$ be positive integers with $r\geq 2$. We say that $r$ families $\mathcal{F}_1\subseteq{V\brack k_1},\mathcal{F}_2\subseteq{V\brack k_2},\ldots,\mathcal{F}_r\subseteq{V\brack k_r}$ are \emph{$r$-cross $t$-intersecting} if $\dim(F_1\cap F_2\cap\cdots\cap F_r)\geq t$	
for all $F_i\in\mathcal{F}_i,\ 1\leq i\leq r.$ For convenience, we just say `cross $t$-intersecting' instead of `$2$-cross $t$-intersecing'.
The $r$-cross $t$-intersecting  families $\mathcal{F}_1$, $\mathcal{F}_2,\ldots,\mathcal{F}_r$ are said to be \emph{trivial} if $\dim(\cap_{1\leq i\leq r}\cap_{F\in\mathcal{F}_i}F)\geq t$ and \emph{non-trivial} otherwise, and said to be \emph{maximal} if $\mathcal{G}_i=\mathcal{F}_i$ $(i=1,2,\ldots,r)$ for any $r$-cross $t$-intersecting families $\mathcal{G}_1\subseteq{V\brack k_1},\mathcal{G}_2\subseteq{V\brack k_2},\ldots,\mathcal{G}_r\subseteq{V\brack k_r}$ with $\mathcal{F}_i\subseteq \mathcal{G}_{i}$ $(i=1,2,\ldots,r)$.

The problem of maximizing the sum or the product of sizes of cross $t$-intersecting families for vector spaces had been well studied so far. Wang and Zhang \cite{Wang-Zhang} completely determined the structure of cross $t$-intersecting families with maximum sum of their sizes. In the following, we focus on the cross $t$-intersecting families with maximum product of their sizes, and let $\mathcal{F}_1\subseteq{V\brack k_1}$ and $\mathcal{F}_2\subseteq{V\brack k_2}$ be two such families. Observe that the structure of these families with $\mathcal{F}_1=\mathcal{F}_2$ can be obtained from the Erd\H{o}s-Ko-Rado Theorem for vector spaces. In \cite{Tokushige-2013}, using the eigenvalue method, Tokushige described the structure of these families with $k_1=k_2$ (see also \cite[Theorem~8.3]{Frankl-Tokushige}). In \cite{Suda-Tanaka}, Suda and Tanaka determined the structure of these families with $t=1$ and $k_1\neq k_2$. Frankl and Tokushige \cite{Frankl-Tokushige} also proposed the following conjecture for $r$-cross $1$-intersecting families.
\begin{conj}\label{con-1}
	Let $(r-1)/r\geq \max\{k_1/n,\ldots,k_r/n\}.$ If $\mathcal{F}_1\subseteq{V\brack k_1},\ldots, \mathcal{F}_r\subseteq{V\brack k_r}$ are $r$-cross $1$-intersecting, then $\prod_{i=1}^r|\mathcal{F}_i|\leq\prod_{i=1}^r{n-1\brack k_i-1}.$
\end{conj}

In this paper, we first study the $r$-cross $t$-intersecting families for vector spaces. The following theorem describes the structure of $r$-cross $t$-intersecting families with maximum product of their sizes, which proves Conjecture~\ref{con-1} with large $n$.
\begin{thm}\label{v-main-1}
	Let $n,r,k_1,k_2,\ldots,k_r$ and $t$ be positive integers with $r\geq 2,$ $k_1\geq k_2\geq\cdots\geq k_r\geq t$ and $n\geq k_1+k_2+t+1$. If $\mathcal{F}_1\subseteq {V\brack k_1}$,  $\mathcal{F}_2\subseteq {V\brack k_2},$\ldots, $\mathcal{F}_r\subseteq {V\brack k_r}$ are $r$-cross $t$-intersecting families, then
	$$
	\prod_{i=1}^r|\mathcal{F}_i|\leq \prod_{i=1}^r{n-t\brack k_i-t}.
	$$
	The equality holds only if $\mathcal{F}_i=\left\{F\in {V\brack k_i}\mid T\subseteq F\right\}\ (i=1,2,\ldots,r)$ for some $T\in {V\brack t}.$
\end{thm}

We also refer the readers to \cite{Borg-2016} for the definition and the results about the $r$-cross $t$-intersecting families for finite sets.

To present our second result let us introduce the following families of subspaces of $V$. Let $k,\ell,s$ be positive integers with $\min\{k,\ell\}\geq s$. Suppose that $X\in{V\brack s}$ and $M\in{V\brack \ell+1}$ with $X\subseteq M$, and $T$ is a subspace of $V$. Write
\begin{align*}
	\mathcal{A}(k,s+1,X,M)=&\left\{F\in{V\brack k}\mid X\subseteq F,\ \dim(F\cap M)\geq s+1\right\},\\
	\mathcal{B}(\ell,X,M)=&\left\{F\in{V\brack \ell}  \mid X\subseteq F \right\}\cup{M\brack \ell},\\
	\mathcal{C}(k,T)=& \left\{F\in{V\brack k}  \mid T\subseteq F \right\},\\
	\mathcal{D}(\ell,s,T)=& \left\{F\in{V\brack \ell}  \mid \dim(F\cap T)\geq s \right\}.
\end{align*}

Let $\mathcal{F}_1\subseteq{V\brack k_1}$ and $\mathcal{F}_2\subseteq{V\brack k_2}$ be non-trivial cross $t$-intersecting families with $k_1\geq k_2$. Suppose that $k_2=t$. Since $\mathcal{F}_1$ and $\mathcal{F}_2$ are nontrivial, we have $k_1>t$ and $|\mathcal{F}_2|\geq 2.$ Let $W$ be the minimal subspace of $V$ containing all elements of $\mathcal{F}_2$. Then $t+1\leq\dim(W)\leq k_1$, $\mathcal{F}_1\subseteq \mathcal{C}(k_1, W)$ and $\mathcal{F}_2\subseteq {W\brack t}.$ Observe that $\mathcal{C}(k_1, W)$ and ${W\brack t}$ are maximal non-trivial cross $t$-intersecting families. Therefore, for any maximal non-trivial cross $t$-intersecting families $\mathcal{G}_1\subseteq{V\brack k_1}$ and $\mathcal{G}_2\subseteq{V\brack k_2}$ with $k_1\geq k_2=t$, there exists a subspace $W$ such that $\mathcal{G}_1=\mathcal{C}(k_1, W)$ and $\mathcal{G}_2={W\brack t}.$ The following theorem describes the structure of $\mathcal{F}_1$ and $\mathcal{F}_2$ with maximum product of their sizes for $k_2\geq t+1.$
\begin{thm}\label{v-main-2}
	Let $n$, $k_1$, $k_2$ and $t$ be positive integers satisfying $k_1\geq k_2\geq t+1$, $n\geq k_1+k_2+t+3$, and $(k_1,k_2,t)\neq(2,2,1),$ $(3,2,1)$ or $(4,2,1)$. Assume that $\mathcal{F}_1\subseteq{V\brack k_1}$ and $\mathcal{F}_2\subseteq{V\brack k_2}$ are non-trivial cross $t$-intersecting families.
	\begin{enumerate}[{\rm (i)}]	
		\item If $k_2\geq 2t+1$, then
		$$
		|\mathcal{F}_1||\mathcal{F}_2|\leq 	\left({n-t\brack k_1-t}-q^{(k_2+1-t)(k_1-t)}{n-k_2-1\brack k_1-t}\right)\left({n-t\brack k_2-t}+q^{k_2+1-t}{t\brack 1}\right),
		$$
		and the equality holds only if there exist a $t$-subspace $X$ and a $(k_2+1)$-subspace $M$ of $V$ with $X\subseteq M$ such that one of the following holds:
		\begin{itemize}
			\item[{\rm(ia)}] $\mathcal{F}_1=\mathcal{A}(k_1,t+1,X,M)$ and $\mathcal{F}_2=\mathcal{B}(k_2, X,M)$,
			\item[{\rm(ib)}] $k_1=k_2$, $\mathcal{F}_1=\mathcal{B}(k_1, X,M)$ and $\mathcal{F}_2=\mathcal{A}(k_2,t+1,X,M)$.
		\end{itemize}
		\item If $t+1\leq k_2\leq 2t$, then
		$$
		|\mathcal{F}_1||\mathcal{F}_2|\leq 	{n-t-1\brack k_1-t-1}\left(q^{k_2-t}{t+1\brack 1}{n-t-1\brack k_2-t}+{n-t-1\brack k_2-t-1}\right),
		$$
		and the equality holds only if there exists a $(t+1)$-subspace $T$ of $V$ such that one of the following holds:
		\begin{itemize}
			\item[{\rm(iia)}] $\mathcal{F}_1=\mathcal{C}(k_1,T)$ and $\mathcal{F}_2=\mathcal{D}(k_2, t, T)$,
			\item[{\rm(iib)}] $k_1=k_2$, $\mathcal{F}_1=\mathcal{D}(k_1, t, T)$ and $\mathcal{F}_2=\mathcal{C}(k_2,T).$
		\end{itemize}	
	\end{enumerate}
\end{thm}

A family $\mathcal{F}\subseteq {V\brack k}$ is called \emph{$r$-wise $t$-intersecting} if for all $F_1,\ldots,F_r\in\mathcal{F}$ one has $\dim(F_1\cap\cdots\cap F_r)\geq t$. Observe that the concept of `$r$-wise $t$-intersecting' is a natural generalization of the classical `$t$-intersecting', and if $\mathcal{F}_1,\mathcal{F}_2,\ldots,\mathcal{F}_r\subseteq{V\brack k}$ are $r$-cross $t$-intersecting with $\mathcal{F}_1=\mathcal{F}_2=\cdots=\mathcal{F}_r=\mathcal{F}$, then $\mathcal{F}$ is $r$-wise $t$-intersecting. Similarly, an $r$-wise $t$-intersecting family $\mathcal{F}\subseteq{V\brack k}$ is called \emph{trivial} if all its members contain a common specified $t$-subspace of $V$ and \emph{non-trivial} otherwise, and is said to be \emph{maximal} if $\mathcal{F}\cup \{F\}$ is not $r$-wise $t$-intersecting for each $F\in{V\brack k}\setminus\mathcal{F}$. In \cite{Chowdhury-2010}, using the shadows in vector space, Chowdhury and Patk\'{o}s determined the structure of extremal $r$-wise $t$-intersecting families of ${V\brack k}$ with $(r-1)n\geq rk$. From the Erd\H{o}s-Ko-Rado Theorem for vector spaces, it is known that if $n>2k$, each extremal $r$-wise $t$-intersecting family of ${V\brack k}$ is trivial. Therefore, in the third main theorem, we will describe extremal non-trivial $r$-wise $t$-intersecting family of ${V\brack k}$ with large $n$.

Let
\begin{align}
h_1(d,k,n)=&{n-d\brack k-d}-q^{(k+1-d)(k-d)}{n-k-1\brack k-d}+q^{k+1-d}{d\brack 1},\label{eqq1}\\
h_2(d,k,n)=&{n-d\brack k-d}-q^{(k-d)^2}{n-k\brack k-d}+q^{k-d+1}{n-k\brack 1}{d\brack 1}.\label{eqq2}
\end{align}
Suppose that $n\geq 2k\geq 6$. From \cite[Lemma~2.6]{Cao-vec}, observe that $h_1(d,k,n)>h_2(d,k,n)$ if $1\leq d<k-2$, and $h_1(d,k,n)<h_2(d,k,n)$ if $d=k-2.$
\begin{thm}\label{v-main-3}
	Let $n$, $k$, $t$ and $r$ be positive integers with $r\geq 3$, $t+r-2\leq k-2$ and $2k+t+r+2\leq n$. Let $\mathcal{F}\subseteq{V\brack k}$ be a non-trivial $r$-wise $t$-intersecting family. Then the following hold.
	\begin{enumerate}[{\rm(i)}]
		\item  Suppose that $t+r-2\leq \frac{k}{2}-1$. Then
		$$
		|\mathcal{F}|\leq {n-t-r+2\brack k-t-r+2}-q^{(k+3-t-r)(k+2-t-r)}{n-k-1\brack k-t-r+2}+q^{k+3-t-r}{t+r-2\brack 1},
		$$
		and equality holds if and only if $\mathcal{F}=\mathcal{A}(k,t+r-1,X,M)\cup{M\brack k}$ for some $(t+r-2)$-subspace $X$ and $(k+1)$-subspace $M$ of $V$ with $X\subset M$. Moreover, if $|\mathcal{F}|>h_2(t+r-2,k,n)$, then $\mathcal{F}$ is a subfamily of $\mathcal{A}(k,t+r-1,X,M)\cup{M\brack k}$ for some $(t+r-2)$-subspace $X$ and $(k+1)$-subspace $M$ of $V$ with $X\subset M$.

		\item Suppose that $\frac{k}{2}-1<t+r-2\leq k-2.$ Then
		$$
		|\mathcal{F}|\leq {t+r\brack 1}{n-t-r+1\brack k-t-r+1}-q{t+r-1\brack 1}{n-t-r\brack k-t-r},
		$$
		and equality holds if and only if $\mathcal{F}=\mathcal{D}(k,t+r-1,Z)$ for some $(t+r)$-subspace $Z$ of $V$. Moreover, if $|\mathcal{F}|>h_2(t+r-2,k,n)$, or $t+r-2=k-2$ and $|\mathcal{F}|>h_1(t+r-2,k,n)$, then $\mathcal{F}$ is a subfamily of $\mathcal{A}(k,t+r-1,X,M)\cup{M\brack k}$  for some $(t+r-2)$-subspace $X$ and $(k+1)$-subspace $M$ of $V$ with $X\subset M$, or a subfamily of  $\mathcal{D}(k,t+r-1,Z)$ for some $(t+r)$-subspace $Z$.
	\end{enumerate}
\end{thm}

In \cite{O-V-2021}, O'Neill and Verstra\"{e}te described the structure of extremal non-trivial $r$-wise $1$-intersecting families of uniform subsets of a set, and also gave a stability theorem. When $t=1$ in Theorem~\ref{v-main-3}, our result can be viewed as a vector space version of the result in \cite{O-V-2021}.
We also refer the readers to \cite{a1,PL,a3,a4} for the upper bounds of an $r$-wise $\mathcal{L}$-intersecting families for vector spaces, where $\mathcal{L}$ is a set of non-negative integers.
\begin{rema}
\label{rem1}
{\em
There is a vast, excellent literature on these intersection problems for finite sets. To limit the scope of this paper, we will not introduce them here. We refer the readers to the monographs \cite{GK} and \cite{FT-2018} for the systematic introduction. In another paper, we will study the non-trivial cross $t$-intersecting families and the non-trivial $r$-wise $t$-intersecting families for finite sets.
}
\end{rema}

The rest of the paper is organized as follows. In Sections 2, 3 and 4, we will prove Theorems~\ref{v-main-1}, \ref{v-main-2} and \ref{v-main-3}, respectively. In Section 5, we will prove some inequalities used in this paper.

\section{$r$-cross $t$-intersecting families with maximum product of their sizes}

In this section, we begin with two useful lemmas which contain some inequalities about the Gaussian coefficients and a formula for counting the number of some special subspaces. Then we describe the structure of the $r$-cross $t$-intersecting families with maximum product of their sizes. The following lemma can be easily proved.
\begin{lemma}\label{lem1-1-1}
	Let $m$ and $i$ be positive integers with $i\leq m.$ Then the following hold.
	\begin{itemize}
		\item[{\rm (i)}] ${m\brack i}={m-1\brack i-1}+q^i{m-1\brack i}$ and ${m\brack i}=\frac{q^m-1}{q^i-1}\cdot{m-1\brack i-1}$.
		\item[{\rm (ii)}] $q^{m-i}<\frac{q^m-1}{q^i-1}<q^{m-i+1}$ and $q^{i-m-1}<\frac{q^i-1}{q^m-1}<q^{i-m}$ if $i <m$.
		\item[{\rm(iii)}] $q^{i(m-i)}\leq{m\brack i}< q^{i(m-i+1)}$, and $q^{i(m-i)}<{m\brack i}$ if $i <m$.
	\end{itemize}
\end{lemma}

Let $W$ be an $(e+\ell)$-dimensional vector space over $\mathbb{F}_q$, where $\ell,e \ge 1$, and let $L$ be a fixed $\ell$-subspace of $W$. We say that an $m$-subspace $U$ is of \emph{type} $(m,h)$ if $\dim(U\cap L)=h$. Define $\mathcal{M}(m,h;e+\ell,e)$ to be the set of all subspaces of $W$ with type $(m,h)$.
Define $N'(m_{1},h_{1};m,h;e+\ell,e)$ to be the number of subspaces of $W$ with type $(m,h)$ containing a given subspace with type $(m_{1},h_{1})$.
\begin{lemma}{\rm{(\cite{KW})}}\label{lem5}
	$N^{'}(m_{1},h_{1};m,h;e+\ell,e)\not= 0$ if and only if $0\leq h_{1}\leq h\leq \ell$ and $0\leq m_{1}-h_{1}\leq m-h\leq e$. Moreover, if $N^{'}(m_{1},h_{1};m,h;e+\ell,e)\neq 0,$ then
	$$ N^{'}(m_{1},h_{1};m,h;e+\ell,e)=q^{(\ell-h)(m-h-m_{1}+h_{1})}{{e-(m_{1}-h_{1})}\brack{(m-h)-(m_{1}-h_{1})}}{{\ell-h_{1}}\brack{h-h_{1}}}.
	$$
\end{lemma}

Observe that $$|\mathcal{M}(m,h;e+\ell,e)|=N'(0,0;m,h;e+\ell,e)=q^{(\ell-h)(m-h)}{e\brack m-h}{\ell\brack h}.$$ Let $a$, $b$ and $c$ be three non-negative integers with $a\leq b\leq c$, and let $A$ and $C$ be two spaces with $\dim(A)=a$, $\dim(C)=c$ and $A\subseteq C$. From \cite[Lemma~9.3.2]{ABE}, we have
\begin{align}\label{eeq1}
\left|\left\{B\in{C\brack b}\mid A\subseteq B\right\}\right|={c-a\brack b-a}.
\end{align}

Set
\begin{align}
	f_1(k,\ell,n,t)=&{t+1\brack 1}{k-t+1\brack 1}{n-t-1\brack \ell-t-1},\label{eqq3}\\
	f_2(m,k,\ell,n,t)=&{m\brack t}{k\brack 1}^{m-t-2}{k-t+1\brack 1}^2{n-m\brack \ell-m}\label{eqq4}.
\end{align}

\begin{lemma} \label{v-upper-L}
	Let $n,\ k,\ \ell,\ t, \ m$ and $s$ be non-negative integers with $t+2\leq m\leq\ell$ and $s+t\leq k$.
	\begin{enumerate}[{\rm(i)}]
		\item If $n\geq k+\ell+2$, we have $f_1(k,\ell,n,t)<{n-t\brack \ell-t}$.
		\item If $n\geq k+\ell+t+1$, the function $f_2(m,k,\ell,n,t)$ is decreasing as $m\in\{t+2, t+3,\ldots,\ell\}$ increases, and $f_2(m,k,\ell,n,t)<f_1(k,\ell,n,t).$
		\item If $n\geq k+\ell$, the function ${\ell-r\brack t-r}{n-s-t+r\brack k-s-t+r}$ is increasing as $r\in\{0,1,\ldots,t-1\}$ increases.
	\end{enumerate}	
\end{lemma}
\proof (i)\quad Since $n\geq k+\ell+2$, observe that
$$
f_1(k,\ell,n,t){n-t\brack \ell-t}^{-1}={t+1\brack 1}{k-t+1\brack 1}\cdot\frac{q^{\ell-t}-1}{q^{n-t}-1}<q^{(t+1)+(k-t+1)+(\ell-n)}\leq 1
$$
from Lemma~\ref{lem1-1-1} (ii) and (iii). Then (i) holds.

(ii)\quad For each $m\in\{t+2,t+3,\ldots,\ell-1\}$, by Lemma~\ref{lem1-1-1} (ii), we have
$$
\frac{f_2(m+1,k,\ell,n,t)}{f_2(m,k,\ell,n,t)}={k\brack 1}\cdot\frac{(q^{m+1}-1)(q^{\ell-m}-1)}{(q^{m-t+1}-1)(q^{n-m}-1)}< q^{k+(t+1)+(\ell-n)}\leq 1.
$$

By Lemma~\ref{lem1-1-1} (ii) and (iii), we have
$$
\frac{f_1(k,\ell,n,t)}{f_2(t+2,k,\ell,n,t)}=\frac{(q^{2}-1)(q-1)(q^{n-t-1}-1)}{(q^{t+2}-1)(q^{k-t+1}-1)(q^{\ell-t-1}-1)}> q^{(-t-1)-(k-t+1)+(n-\ell)}\geq 1,
$$
implying that (ii) holds.

(iii)\quad Let
$$
g(r)={\ell-r\brack t-r}{n-s-t+r\brack k-s-t+r}
$$
for $r\in\{0,1,2,\ldots,t-1\}$. Since $n\geq k+\ell$, by Lemma~\ref{lem1-1-1} (ii), we have
$$
\frac{g(r+1)}{g(r)}=\frac{(q^{t-r}-1)(q^{n-s-t+r+1}-1)}{(q^{\ell-r}-1)(q^{k-s-t+r+1}-1)}>q^{(t-\ell-1)+(n-k)}>1
$$
for each $r\in\{0,1,\ldots,t-2\}$. That is, the function $g(r)$ is increasing as $r\in\{0,1,\ldots,t-1\}$ increases.    $\qed$

For two subspaces $A$ and $B$, recall that the set $\{a+b\in V\mid a\in A,\ b\in B\}$ is a subspace of $V$, and is usually denoted by $A+B$. For a family $\mathcal{F}\subseteq{V\brack k}$ and a subspace $S$ of $V$, let $\mathcal{F}_S$ denote the set of $k$-subspaces in $\mathcal{F}$ which contain $S$.
\begin{lemma}\label{v-prer}
Let $n,\ k,\ \ell,\ t$ and $s$ be non-negative integers with $n\geq k+\ell$, and $G\in{V\brack \ell}$ and $S\in{V\brack s}$ be two subspaces of $V$ with $\dim(G\cap S)=r<t$. Let $\mathcal{F}\subseteq{V\brack k}$ be a family satisfying $\dim(G\cap F)\geq t$ for all $F\in\mathcal{F}$. Then for each $i\in\{1,2,\ldots,t-r\}$, there exists an $(s+i)$-subspace $U_i$ with $S\subseteq U_i$  such that  $|\mathcal{F}_S|\leq {\ell-r\brack i} |\mathcal{F}_{U_i}|$. Furthermore,  we have $|\mathcal{F}_{S}|\leq {\ell-r\brack t-r}{n-s-t+r\brack k-s-t+r}$.
\end{lemma}
\proof If $\mathcal{F}_S=\emptyset$, the required result is clear. Now suppose that $\mathcal{F}_S\neq \emptyset$. For each $i\in\{1,2,\ldots, t-r\}$, write
$$
\mathcal{H}_i=\left\{H\in{G+S\brack s+i}\mid S\subseteq H\right\}.
$$
Observe that $\dim(G+S)=\ell+s-r$. For each $F\in\mathcal{F}_S$, since $\dim(F\cap G)\geq t$, we have $\dim(F+S+G)=\dim(F+G)=\dim(F)+\dim(G)-\dim(F\cap G)\leq k+\ell-t,$ implying that $\dim(F\cap (S+G))=\dim(F)+\dim(S+G)-\dim(F+S+G)\geq s+t-r.$ It follows that there exists $H\in\mathcal{H}_i$ such that $H\subseteq F.$ Therefore, $\mathcal{F}_S=\cup_{H\in\mathcal{H}_i}\mathcal{F}_{H}$. Let $U_i$ be an element in $\mathcal{H}_{i}$ such that $|\mathcal{F}_{U_i}|=\max\{|\mathcal{F}_{H}|\mid  H\in\mathcal{H}_i\}$. By (\ref{eeq1}), we have $|\mathcal{H}_i|={\ell-r\brack i}$, and then the former part holds. Setting $i=t-r$, then $\dim(U_i)=s+t-r$, and the latter part holds due to $|\mathcal{F}_{U_{t-r}}|\leq {n-s-t+r\brack k-s-t+r}$.   $\qed$

Let $\mathcal{F}\subseteq {V\brack k}$ and $\mathcal{G}\subseteq {V\brack \ell}$ be cross $t$-intersecting families. Obviously, any element in $\mathcal{G}$ is a $t$-cover of $\mathcal{F}$, and any element in $\mathcal{F}$ is a $t$-cover of $\mathcal{G}$. Then $t\leq\tau_{t}(\mathcal{F})\leq \ell$ and $t\leq\tau_{t}(\mathcal{G})\leq k$.
\begin{lemma}\label{v-upper-F}
Let $n,$ $k,$ $\ell$ and $t$ be non-negative integers with $n\geq k+\ell$, and $\mathcal{F}\subseteq {V\brack k}$ and $\mathcal{G}\subseteq {V\brack \ell}$ be maximal cross $t$-intersecting families. Assume that $\tau_{t}(\mathcal{F})=m_k$  and $\tau_{t}(\mathcal{G})=m_\ell$. Then
	$$
	|\mathcal{F}|\leq{m_k\brack t}{n-t\brack k-t}.
	$$
Moreover, the following hold.
\begin{enumerate}[{\rm(i)}]
	\item If $m_\ell=t+1$, then
	$$
	|\mathcal{F}|\leq{m_k\brack t}{\ell-t+1\brack 1}{n-t-1\brack k-t-1}.
	$$
	\item If $m_\ell\geq t+2$, then
	$$
	|\mathcal{F}|\leq {m_k\brack t}{\ell\brack 1}^{m_\ell-t-2}{\ell-t+1\brack 1}^2{n-m_\ell\brack k-m_\ell}.
	$$
\end{enumerate}
\end{lemma}
\proof Let $T$ be a $t$-cover of $\mathcal{F}$ with $\dim(T)=m_k.$ Then for each $F\in\mathcal{F},$ we have $\dim(F\cap T)\geq t$ and there exists $H\in{T\brack t}$ such that $H\subseteq F,$ implying that $\mathcal{F}\subseteq \cup_{H\in{T\brack t}}\mathcal{F}_{H}$. Let  $H_1$ be a $t$-subspace of $T$ such that $|\mathcal{F}_{H_1}|=\max\{|\mathcal{F}_{H}|\mid H\in{T\brack t}\}$. We have
\begin{align}\label{v-upp-1}
	|\mathcal{F}|\leq{m_k\brack t}|\mathcal{F}_{H_1}|.
\end{align}
Since $|\mathcal{F}_{H_1}|\leq {n-t\brack k-t},$ we have $|\mathcal{F}|\leq{m_k\brack t}{n-t\brack k-t}$ by (\ref{v-upp-1}).

\

(i)\quad Suppose that $m_\ell=t+1$. It follows from $\dim(H_1)=t$ that there exists $G_1\in\mathcal{G}$ such that $\dim(G_1\cap H_1)=r<t$. By Lemmas~\ref{v-prer} and \ref{v-upper-L} (iii), we have $$|\mathcal{F}_{H_1}|\leq{\ell-r\brack t-r}{n-2t+r\brack k-2t+r}\leq{\ell-t+1\brack 1}{n-t-1\brack k-t-1},$$ implying that (i) holds from (\ref{v-upp-1}).

\

(ii)\quad Suppose that $m_{\ell}\geq t+2$. Firstly we claim that there exists $H^\prime\in{V\brack m_\ell-2}$ such that
\begin{align}\label{v-upp-2}
|\mathcal{F}|\leq{m_k\brack t}{\ell\brack 1}^{m_\ell-t-2}|\mathcal{F}_{H^\prime}|.	
\end{align}
If $m_\ell=t+2$, it is clear that (\ref{v-upp-2}) holds by setting $H^\prime=H_1$. If $m_\ell\geq t+3$, using Lemma~\ref{v-prer} repeatedly, then there exist $$H_2\in{V\brack t+1},\ H_3\in{V\brack t+2},\ldots,H_{m_\ell-t-1}\in{V\brack m_\ell-2}$$ such that $H_i\subseteq H_{i+1}$ and $|\mathcal{F}_{H_{i}}|\leq {\ell\brack 1}|\mathcal{F}_{H_{i+1}}|$ for each $i\in\{1,2,\ldots,m_\ell-t-2\}$, which implies that (\ref*{v-upp-2}) holds by setting $H^\prime=H_{m_\ell-t-1}$. Therefore, the claim holds.

Since $\tau_t(\mathcal{G})>m_\ell-2,$ there exists $G_2\in\mathcal{G}$ such that $\dim(H^\prime\cap G_2)<t.$ If $\dim(H^\prime\cap G_2)=r\leq t-2$, by Lemmas~\ref{v-prer} and \ref{v-upper-L} (iii), we have
\begin{align}\label{v-upp-3}
|\mathcal{F}_{H^\prime}|\leq {\ell-r\brack t-r}{n-m_\ell+2-t+r\brack k-m_\ell+2-t+r}\leq {\ell-t+2\brack 2}{n-m_\ell\brack k-m_\ell}.
\end{align}
Suppose that $\dim(H^\prime\cap G_2)=t-1$. By Lemma~\ref{v-prer}, there exists an $(m_\ell-1)$-subspace $H^{\prime\prime}$ such that $|\mathcal{F}_{H^\prime}|\leq{\ell-t+1\brack 1}|\mathcal{F}_{H^{\prime\prime}}|$. Since $\tau_t(\mathcal{G})>m_\ell-1$,  there exists $G_3\in\mathcal{G}$ such that $\dim(H^\prime\cap G_3)<t.$ If $\dim(H^\prime\cap G_3)\leq t-2$, by Lemma~\ref{v-prer} again, we have
\begin{align}\label{v-upp-4}
|\mathcal{F}_{H^\prime}|\leq {\ell-t+1\brack 1}|\mathcal{F}_{H^{\prime\prime}}|\leq {\ell-t+1\brack 1}{\ell-t+2\brack 2}{n-m_\ell-1\brack k-m_\ell-1}.
\end{align}
If $\dim(H^\prime\cap G_3)=t-1$, by Lemma~\ref{v-prer}, then
\begin{align}\label{v-upp-5}
|\mathcal{F}_{H^\prime}|\leq {\ell-t+1\brack 1}|\mathcal{F}_{H^{\prime\prime}}|\leq {\ell-t+1\brack 1}^2{n-m_\ell\brack k-m_\ell}.
\end{align}

By Lemma~\ref{lem1-1-1} (ii) and $n\geq k+\ell,$  we have $q^{\ell-t+1}>\frac{q^{\ell-t+2}-1}{q+1}$ and
$$
\frac{q^{n-m_{\ell}}-1}{q^{k-m_{\ell}}-1}>q^{n-k}\geq q^{\ell}>q^{\ell-t+1}>\frac{q^{\ell-t+2}-1}{q^2-1},
$$
implying that
$$
{\ell-t+1\brack 1}^2{n-m_\ell\brack k-m_\ell}\geq\max\left\{{\ell-t+2\brack 2}{n-m_\ell\brack k-m_\ell},\ {\ell-t+1\brack 1}{\ell-t+2\brack 2}{n-m_\ell-1\brack k-m_\ell-1}\right\}.
$$
This together with (\ref{v-upp-2}), (\ref{v-upp-3}), (\ref{v-upp-4}) and (\ref{v-upp-5}) yields (ii) holds. $\qed$
\begin{cor}\label{v-cor}
Let $n$, $k$, $\ell$ and $t$ be non-negative integers with $n\geq k+\ell+t+1$, and $\mathcal{F}\subseteq{V\brack k}$ and $\mathcal{G}\subseteq {V\brack \ell}$ be maximal cross $t$-intersecting families. Assume that $\tau_t(\mathcal{F})=m_k$ and $\tau_t(\mathcal{G})=m_\ell$ with $m_k\leq m_\ell$.
\begin{enumerate}[{\rm (i)}]
\item If $m_k=t,$ then
$$
|\mathcal{F}||\mathcal{G}|\leq\left\{
\begin{array}{ll}
{n-t\brack \ell-t}f_1(\ell,k,n,t), & \mbox{if}\ m_\ell=t+1,\\
{n-t\brack \ell-t}f_2(m_\ell,\ell,k,n,t), & \mbox{if}\ m_\ell\geq t+2.
\end{array}
\right.	
$$
\item If $m_k\geq t+1,$ then $|\mathcal{F}||\mathcal{G}|<f_1(k,\ell,n,t)f_1(\ell,k,n,t)$.
\end{enumerate}
\end{cor}
\proof Applying Lemma~\ref{v-upper-F} to $\mathcal{F}$ and $\mathcal{G}$ respectively,  we have
\begin{align}\label{v-upper-F1}
	|\mathcal{F}|\leq\left\{
	\begin{array}{ll}
		{m_k\brack t}{\ell-t+1\brack 1}{n-t-1\brack k-t-1},&\mbox{if}\ m_\ell=t+1,\vspace{0.1cm}\\
		{m_k\brack t}{\ell\brack 1}^{m_\ell-t-2}{\ell-t+1\brack 1}^2{n-m_\ell\brack k-m_\ell},&\mbox{if}\ m_\ell\geq t+2,
	\end{array}
	\right.
\end{align}
and
\begin{align}\label{v-upper-F2}
	|\mathcal{G}|\leq\left\{
	\begin{array}{ll}
		{m_\ell\brack t}{n-t\brack \ell-t},&\mbox{if}\ m_k=t,\vspace{0.1cm}\\
		{m_\ell\brack t}{k-t+1\brack 1}{n-t-1\brack \ell-t-1},&\mbox{if}\ m_k=t+1,\vspace{0.1cm}\\
		{m_\ell\brack t}{k\brack 1}^{m_k-t-2}{k-t+1\brack 1}^2{n-m_k\brack \ell-m_k},&\mbox{if}\ m_k\geq t+2.
	\end{array}
	\right.
\end{align}

(i)\quad It is straightforward to verify that (i) holds from (\ref{v-upper-F1}) and (\ref{v-upper-F2}).

(ii)\quad By (\ref{v-upper-F1}) and (\ref{v-upper-F2}) again, we have

$$
|\mathcal{F}||\mathcal{G}|\leq \left\{
\begin{array}{ll}
	f_1(k,\ell,n,t)f_1(\ell,k,n,t),&\mbox{if}\ m_k=t+1,\ m_\ell=t+1,\vspace{0.1cm}\\
	f_1(k,\ell,n,t)f_2(m_\ell,\ell,k,n,t),&\mbox{if}\ m_k=t+1,\ m_\ell\geq t+2,\vspace{0.1cm}\\
	f_2(m_k,k,\ell,n,t)f_2(m_\ell,\ell,k,n,t),&\mbox{if}\ m_\ell\geq m_k\geq t+2.
\end{array}
\right.
$$
From Lemma~\ref{v-upper-L} (ii), we have that $|\mathcal{F}||\mathcal{G}|<f_1(k,\ell,n,t)f_1(\ell,k,n,t)$ holds.  $\qed$
\begin{lemma}\label{v-FTt}
	Let $n,k_1,k_2$ and $t$ be positive integers with $n\geq k_1+k_2$, and $\mathcal{F}_1\subseteq {V\brack k_1}$ and $\mathcal{F}_2\subseteq {V\brack k_2}$ be maximal cross $t$-intersecting families. For $i=1,2$, suppose that $\mathcal{T}_i$ is the set of the $t$-covers of $\mathcal{F}_i$ with dimensional $\tau_t(\mathcal{F}_i)$. Then $\mathcal{T}_1$ and $\mathcal{T}_2$ are cross $t$-intersecting families.
\end{lemma}
\proof It follows from the maximality of $\mathcal{F}_1$ and $\mathcal{F}_2$ that,
for each $T_1\in\mathcal{T}_1$, $\mathcal{F}_2$ contains all $k_2$-subspaces of $V$ containing $T_1;$ and for each $T_2\in\mathcal{T}_2$, $\mathcal{F}_1$ contains all $k_1$-subspaces of $V$ containing $T_2.$ If there exist $A_1\in\mathcal{T}_1$ and $A_2\in\mathcal{T}_2$ such that $\dim(A_1\cap A_2)<t,$ by $n\geq k_1+k_2$, then there must exist $F_1\in\mathcal{F}_1$ and $F_2\in\mathcal{F}_2$ such that $A_1\subseteq F_1$, $A_2\subseteq F_2$ and $\dim(F_1\cap F_2)<t.$ This is impossible because $\mathcal{F}_1$ and $\mathcal{F}_2$ are cross $t$-intersecting families. Hence, $\mathcal{T}_1$ and $\mathcal{T}_2$ are cross $t$-intersecting families.  $\qed$

\noindent{\em \textbf{Proof of Theorem~\ref{v-main-1}.}}\quad We first prove that the theorem holds for $r=2$. It is clear that $|\mathcal{F}_1||\mathcal{F}_2|\leq{n-t\brack k_1-t}{n-t\brack k_2-t}$ if $\mathcal{F}_1$ and $\mathcal{F}_2$ are trivial cross $t$-intersecting families, and equality holds only if  $\mathcal{F}_i=\{F\in{V\brack k_i}\mid T\subseteq F\}\ (i=1,2)$ for some $t$-subspace $T$ of $V.$ Then, it suffices to prove that $|\mathcal{F}_1||\mathcal{F}_2|<{n-t\brack k_1-t}{n-t\brack k_2-t}$ if $\mathcal{F}_1$ and $\mathcal{F}_2$ are non-trivial cross $t$-intersecting families. Without loss of generality, assume that $\mathcal{F}_1$ and $\mathcal{F}_2$ are maximal non-trivial cross $t$-intersecting families with $\tau_t(\mathcal{F}_1)\leq \tau_t(\mathcal{F}_2)$. Set $\tau_t(\mathcal{F}_i)=m_i$ for each $i\in\{1,2\}$. If $m_2=t$, then $m_1=t$, and $\mathcal{T}_1\cap\mathcal{T}_2\neq\emptyset$ by Lemma~\ref{v-FTt}, implying that there exists a $t$-subspace $T$ of $V$ such that $T$ is contained in each member of $\mathcal{F}_1\cup\mathcal{F}_2$. That is impossible because $\mathcal{F}_1$ and $\mathcal{F}_2$ are non-trivial cross $t$-intersecting. If $m_2\geq t+1,$  by Corollary~\ref{v-cor} and Lemma~\ref{v-upper-L} (i) and (ii), we obtain $|\mathcal{F}_1||\mathcal{F}_2|<{n-t\brack k_1-t}{n-t\brack k_2-t}$ as required.

Now suppose that $r\geq 3$. For $u,v\in\{1,2,\ldots,r\}$ with $u<v$, it follows from $\mathcal{F}_u$ and $\mathcal{F}_v$ are cross $t$-intersecting family that
$$
|\mathcal{F}_u||\mathcal{F}_v|\leq{n-t\brack k_u-t}{n-t\brack k_v-t}.
$$
Then we have
{\small\begin{align}\label{v-equ-11}
\left(\prod_{i=1}^r|\mathcal{F}_i|\right)^{r-1}=\prod_{1\leq u<v\leq r}|\mathcal{F}_u||\mathcal{F}_v|\leq\prod_{1\leq u<v\leq r}{n-t\brack k_u-t}{n-t\brack k_v-t}=\left(\prod_{i=1}^r{n-t\brack k_i-t}\right)^{r-1},
\end{align}}
which implies that
$$
\prod_{i=1}^r|\mathcal{F}_i|\leq \prod_{i=1}^r{n-t\brack k_i-t}.
$$

Suppose that $\mathcal{F}_1,\mathcal{F}_2,\ldots,\mathcal{F}_r$ are $r$-cross $t$-intersecting families satisfying $\prod_{i=1}^r|\mathcal{F}_i|= \prod_{i=1}^r{n-t\brack k_i-t}.$ For each $u,v\in\{1,2,\ldots,r\}$ with $u<v$, note that $|\mathcal{F}_u||\mathcal{F}_v|={n-t\brack k_u-t}{n-t\brack k_v-t}$ by (\ref{v-equ-11}), and so $\mathcal{F}_h=\{F\in{V\brack k_h}\mid T_{u,v}\subseteq F\}$ ($h\in\{u,v\}$) for some $t$-subspace $T_{u,v}$ of $V$ due to the theorem holds for $r=2$. From the structure of $\mathcal{F}_i$ $(i\in\{1,2,\ldots,r\})$, observe that there exists a $t$-subspace $T$ of $V$ such that $T_{u,v}=T$ for all $u,v\in\{1,2,\ldots,r\}$ with $u<v$, hence the desired result follows.     $\qed$

\section{Non-trivial cross $t$-intersecting families with maximum product of their sizes}
In this section, we describe the structure of the non-trivial $r$-cross $t$-intersecting families with maximum product of their sizes. Write
{\small\begin{align}
	g_1(k,\ell,n,t)=&{n-t\brack k-t}-q^{(\ell+1-t)(k-t)}{n-\ell-1\brack k-t},\label{eqq5}\\
	g_2(\ell,n,t)=&{n-t\brack \ell-t}+q^{\ell+1-t}{t\brack 1},\label{eqq6}\\
	g_3(k,\ell,n,t)=&{n-t-1\brack k-t-1}\left(q^{\ell-t}{t+1\brack 1}{n-t-1\brack \ell-t}+{n-t-1\brack \ell-t-1}\right),\label{eqq8}\\
	 f_3(k,\ell,n,t)=&{\ell-t\brack 1}{n-t-1\brack k-t-1}+q^{2(\ell-t)}{n-t-2\brack k-t-2},\label{eqq10}\\
	f_4(k,\ell,n,t)=&{n-t\brack \ell-t}+q^{\ell-t+1}{t\brack 1}{k-t\brack 1}{n-t-2\brack \ell-t-1}.\label{eqq11}
\end{align}}
\begin{lemma}\label{lem3-1}
Let $n,\ k,\ \ell$ and $t$ be non-negative integers with $n\geq k+\ell+t+2$, and $\mathcal{F}\subseteq {V\brack k}$ and $\mathcal{G}\subseteq {V\brack \ell}$ be maximal cross $t$-intersecting families with $\tau_t(\mathcal{F})=t$ and $\tau_t(\mathcal{G})=t+1$. Let $X$ be a $t$-cover of $\mathcal{F}$ with $\dim(X)=t$, and $T$ a $t$-cover of $\mathcal{G}$ with $\dim(T)=t+1$.
\begin{enumerate}[{\rm(i)}]
\item For each $F\in\mathcal{F}$ and $G\in\mathcal{G}\setminus\mathcal{G}_{X}$, we have $X\subseteq F$, $\dim(G\cap X)=t-1$ and $\dim(F\cap(G+X))\geq t+1.$
\item If there exists $(\ell+1)$-subspace $M$ of $V$ such that $G+X=M$ for each $G\in\mathcal{G}\setminus\mathcal{G}_{X}$, then
\begin{align*}
	 \mathcal{F}=\mathcal{A}(k,t+1,X,M)\quad\mbox{and}\quad\mathcal{G}=\mathcal{B}(\ell,X,M).
\end{align*}
\item If $T\subseteq F$ for each $F\in\mathcal{F}$, then
\begin{align*}
	\mathcal{F}=\mathcal{C}(k,T)\quad\mbox{and}\quad
	\mathcal{G}=\mathcal{D}(\ell,t,T).
\end{align*}
\item If there exist $G_1,G_2\in\mathcal{G}\setminus\mathcal{G}_X$ such that $G_1+X\neq G_2+X$, and there exists $F^\prime\in\mathcal{F}$ such that $T\nsubseteq F^{\prime}$, then
\begin{align*}
	|\mathcal{F}||\mathcal{G}|\leq f_3(k,\ell,n,t)f_4(k,\ell,n,t).
\end{align*}	
\end{enumerate}
\end{lemma}
\proof (i)\quad Since $X$ is a $t$-cover of $\mathcal{F}$ with $\dim(X)=t$, it is clear that $X\subseteq F$ for each $F\in\mathcal{F}.$ By Lemma~\ref{v-FTt}, observe that $X\subseteq T$ and $\dim(T)-\dim(X)=1.$ For each $G\in\mathcal{G}\setminus\mathcal{G}_{X}$, since $\dim(G\cap T)\geq t$, $X\subseteq T$ and $X\nsubseteq G,$ we have $t=\dim(T\cap G)\leq\dim(X\cap G)+1$ and $\dim(X\cap G)\leq t-1$, implying that $\dim(X\cap G)=t-1$. It follows from $(F\cap G)+(F\cap X)\subseteq F\cap(G+X)$ that
$$
\dim(F\cap(G+X))\geq\dim(F\cap G)+\dim(F\cap X)-\dim(F\cap G\cap X)\geq t+1.
$$

(ii)\quad It is routine to check that $\mathcal{A}(k,t+1,X,M)$ and $\mathcal{B}(\ell,X,M)$ are cross $t$-intersecting with $\tau_t(\mathcal{A}(k,t+1,X,M))=t$ and $\tau_t(\mathcal{B}(\ell,X,M))=t+1.$ For each $F\in\mathcal{F},$ by (i), $F\in\mathcal{A}(k,t+1,X,M)$ holds. For each $G\in\mathcal{G}$, if $X\subseteq G$, then $G\in\mathcal{B}(\ell,X,M);$ if $X\nsubseteq G,$ then $G+X=M$, implying that $G\in\mathcal{B}(\ell,X,M)$. By the maximality, we have that $\mathcal{F}=\mathcal{A}(k,t+1,X,M)$ and $\mathcal{G}=\mathcal{B}(\ell,X,M).$

(iii)\quad It is routine to check that $\mathcal{F}\subseteq \mathcal{C}(k,T)$, $\mathcal{G}\subseteq\mathcal{D}(\ell,t,T)$, and $\mathcal{C}(k,T)$ and $\mathcal{D}(\ell,t,T)$ are cross $t$-intersecting. Hence, we have that (iii) holds.

(iv)\quad Assume that $G_1+X=M_1$, $G_2+X=M_2$, $M^\prime=M_1\cap M_2$ and $\dim(M^\prime)=m.$ Observe that $t\leq m\leq\ell.$ Since $X\subseteq F$ for each $F\in\mathcal{F}$, we have $F^\prime\cap T=X$. Set
\begin{align*}
	\mathcal{F}^\prime=&\left\{F\in{V\brack k}\mid X\subseteq F,\ \dim(F\cap M_1)\geq t+1,\ \dim(F\cap M_2)\geq t+1\right\}, \\
	\mathcal{G}^\prime=&\left\{G\in{V\brack \ell}\mid X\nsubseteq G,\ \dim(G\cap T)=t,\ \dim(G\cap F^\prime)\geq t\right\}.
\end{align*}
By (i) and the assumption, we have that $\mathcal{F}\subseteq\mathcal{F}^\prime$ and $\mathcal{G}\subseteq\{G\in{V\brack \ell}\mid X\subseteq G\}\cup \mathcal{G}^\prime.$ It is clear that $|\{G\in{V\brack \ell}\mid X\subseteq G\}|={n-t\brack \ell-t}$. In order to prove that (iv) holds, it suffices to give the upper bounds of $|\mathcal{F}^\prime|$ and $|\mathcal{G}^\prime|$.

\;

\noindent{\textbf{Step 1.}}\quad Show that $|\mathcal{F}^\prime|\leq f_3(k,\ell,n,t)$.

Set
$$
\mathcal{W}=\left\{(W_1,W_2)\in{M_1\brack t+1}\times{M_2\brack t+1}\mid X\subseteq W_1\nsubseteq M^\prime,\ X\subseteq W_2\nsubseteq M^\prime \right\}.
$$
If there exists $(W_1,W_2)\in\mathcal{W}$ such that $W_1=W_2,$ then $W_1=W_2\subseteq M_1\cap M_2=M^\prime$, a contradiction. Hence, $\dim(W_1\cap W_2)=t$ and $\dim(W_1+W_2)=t+2$ for each $(W_1,W_2)\in\mathcal{W}.$ For each $F\in\mathcal{F}^\prime,$ if $\dim(F\cap M^\prime)\geq t+1,$ then there exists $H\in{M^\prime \brack t+1}$ such that $H\subseteq F;$ if $\dim(F\cap M^\prime)=t,$ from $\dim(F\cap M_1)\geq t+1$ and $\dim(F\cap M_2)\geq t+1$, then there exist $(W_1,W_2)\in\mathcal{W}$ such that $W_1\subseteq F$ and $W_2\subseteq F.$ Therefore, we have
$$
\mathcal{F}^\prime\subseteq \left(\bigcup_{H\in{M^\prime\brack t+1},\ X\subseteq H}\mathcal{F}^\prime_H\right)\bigcup\left(\bigcup_{(W_1,W_2)\in\mathcal{W}}\mathcal{F}^\prime_{W_1+W_2}\right).
$$
Observe that $|\mathcal{F}^\prime_{H}|\leq{n-t-1\brack k-t-1}$ for each $H\in{M^\prime \brack t+1}$ with $X\subseteq H$, $|\mathcal{W}|=\left({\ell+1-t\brack 1}-{m-t\brack 1}\right)^2$ and $|\mathcal{F}^\prime_{W_1+W_2}|\leq{n-t-2\brack k-t-2}$ for each $(W_1,W_2)\in\mathcal{W}$. For each $m\in\{t,t+1,\ldots,\ell\}$, write
$$
f^{\prime}(n,k,\ell,m,t)={m-t\brack 1}{n-t-1\brack k-t-1}+\left({\ell+1-t\brack 1}-{m-t\brack 1}\right)^2{n-t-2\brack k-t-2}.
$$
Then $|\mathcal{F}^\prime|\leq f^{\prime}(n,k,\ell,m,t).$ Observe that
\begin{align*}
f^{\prime}(n,k,\ell,m,t)=&{m-t\brack 1}\left({n-t-1\brack k-t-1}-2{\ell+1-t\brack 1}{n-t-2\brack k-t-2}\right)\\&+\left({\ell+1-t\brack 1}^2+{m-t\brack 1}^2\right){n-t-2\brack k-t-2}.
\end{align*}
Since $n\geq k+\ell+t+1,$ by Lemma~\ref{lem1-1-1} (ii), we have
$$
{n-t-1\brack k-t-1}-2{\ell+1-t\brack 1}{n-t-2\brack k-t-2}\geq q^{n-k}{n-t-2\brack k-t-2}-2q^{\ell+1-t}{n-t-2\brack k-t-2}>0
$$
implying that $f^{\prime}(n,k,\ell,m,t)$ increase as $m\in\{t,t+1,\ldots,\ell\}$ increases. Hence $|\mathcal{F}^\prime|\leq f^{\prime}(n,k,\ell,\ell,t)=f_3(k,\ell,n,t)$ as required.

\;

\noindent{\textbf{Step 2.}}\quad Show that $|\mathcal{G}^\prime|\leq q^{\ell-t+1}{t\brack 1}{k-t\brack 1}{n-t-2\brack \ell-t-1}$.

Set
$$
\mathcal{W}^\prime=\left\{W\in{T+F^\prime\brack t+1}\mid \dim(W\cap T)=t,\ W\cap T\neq X \right\}.
$$
For $G\in\mathcal{G}^\prime,$ since $X\nsubseteq G,$ $\dim(G\cap T)=t$ and $\dim(G\cap F^\prime)\geq t,$ we have $X\nsubseteq G\cap(T+F^\prime)$, and
$$
\dim(G\cap(T+F^\prime))\geq \dim((G\cap T)+(G\cap F^\prime))= \dim(G\cap T)+\dim(G\cap F^\prime)-\dim(G\cap T\cap F^\prime)\geq t+1,
$$
due to $\dim(G\cap T\cap F^\prime)=\dim(G\cap X)\leq t-1$, implying that there exist $W\in\mathcal{W}^\prime$ such that $W\subseteq G.$
Therefore,
$$
\mathcal{G}^\prime\subseteq\bigcup_{W\in\mathcal{W}^\prime}\mathcal{G}^\prime_{W}\subseteq \bigcup_{W\in\mathcal{W}^\prime}\left\{G\in{V\brack \ell}\mid W\subseteq G,\ T\nsubseteq G\right\}.
$$
Observe that
$$
\mathcal{W}^\prime=\left\{W\in{T+F^\prime\brack t+1}\mid \dim(W\cap T)=t\right\}\setminus\left\{W\in{T+F^\prime\brack t+1}\mid W\cap T=X \right\},
$$
and so
$$
|\mathcal{W}^\prime|=N^\prime(0,0;t+1,t;k+1,k-t)- N^\prime(t,t;t+1,t;k+1,k-t)=q^2{t\brack 1}{k-t\brack 1}.
$$

For each $W\in\mathcal{W}^\prime$, since
$$
\left\{G\in{V\brack \ell}\mid W\subseteq G,\ T\nsubseteq G\right\}=\left\{G\in{V\brack \ell}\mid W\subseteq G\right\}\setminus\left\{G\in{V\brack \ell}\mid W\subseteq G,\ T\subseteq G\right\},
$$
we have
$$
\left|\left\{G\in{V\brack \ell}\mid W\subseteq G,\ T\nsubseteq G\right\}\right|={n-t-1\brack\ell-t-1}-{n-t-2\brack \ell-t-2}=q^{\ell-t-1}{n-t-2\brack \ell-t-1}
$$
 Hence, we have that
$$
|\mathcal{G}^\prime|\leq q^{\ell-t+1}{t\brack 1}{k-t\brack 1}{n-t-2\brack \ell-t-1}
$$
as desired.     $\qed$

\begin{lemma}\label{lem-abcd}
Let $k,$ $\ell$ and $t$ be positive integers with $\min\{k,\ell\}\geq t+1$, and $X\in{V\brack t},\ M\in{V\brack \ell+1},\ T\in{V\brack t+1}$ with $X\subseteq M.$	 Then the following hold.
\begin{enumerate}[{\rm(i)}]
\item $|\mathcal{A}(k,t+1,X,M)|=g_1(k,\ell,n,t)$ and $|\mathcal{B}(\ell,X,M)|=g_2(\ell,n,t)$.
\item $|\mathcal{C}(k,T)|={n-t-1\brack k-t-1}$ and $|\mathcal{D}(\ell,t,T)|=q^{\ell-t}{t+1\brack 1}{n-t-1\brack \ell-t}+{n-t-1\brack \ell-t-1}$. Moreover, we have $|\mathcal{C}(k,T)||\mathcal{D}(\ell,t,T)|=g_3(k,\ell,n,t)$.	
\end{enumerate}	
\end{lemma}
\proof (i) Since $$\mathcal{A}(k,t+1,X,M)=\left\{F\in{V\brack k}\mid X\subseteq F\right\}\setminus\left\{F\in{V\brack k}\mid X\subseteq F,\ \dim(F\cap M)= t\right\},$$
we have $|\mathcal{A}(k,t+1,X,M)|={n-t\brack k-t}-N^\prime(t,t;k,t;n,n-\ell-1)$ by (\ref{eeq1})
and the required result follows by Lemma~\ref{lem5}. Observe that
$$
\left|{M\brack\ell}\setminus \left\{F\in{V\brack \ell}\mid X\subseteq F\right\}\right|={\ell+1\brack \ell}-{\ell+1-t\brack\ell-t}=q^{\ell+1-t}{t\brack 1},
 $$and then $|\mathcal{B}(\ell,X,M)|=g_2(\ell,n,t)$ holds from $\left|\left\{F\in{V\brack \ell}\mid X\subseteq F\right\}\right|={n-t\brack \ell-t}$.

(ii) It is clear that $|\mathcal{C}(k,T)|={n-t-1\brack k-t-1}$ holds from the construction of $\mathcal{C}(k,T)$ and (\ref{eeq1}). Since
$$
\mathcal{D}(\ell,t,T)=\left\{F\in{V\brack \ell} \mid \dim(F\cap T)=t \right\}\cup \left\{F\in{V\brack \ell} \mid T\subseteq F \right\},
$$
we have $|\mathcal{D}(\ell,t,T)|=N^\prime(0,0;\ell,t;n,n-t-1)+{n-t-1\brack \ell-t-1},$
and the required result follows by Lemma~\ref{lem5}.  $\qed$

\
\medskip
\

\noindent{\em\textbf{Proof of Theorem~\ref{v-main-2}.}}\quad By Lemmas~\ref{lem-upp-1}, \ref{v-inequ-2}, \ref{v-inequ-3} and  \ref{v-inequ-3-1}, we first have
\begin{align}\label{v-ineq-10}
g_1(k_1,k_2,n,t)g_2(k_2,n,t)\left\{\begin{array}{ll}
>g_3(k_1,k_2,n,t),& \mbox{if}\ k_2\geq 2t+1,\\
<g_3(k_1,k_2,n,t),& \mbox{if}\  k_2\leq 2t.
\end{array}\right.
\end{align}
From Lemma~\ref{lem-abcd}, observe that $|\mathcal{F}_1||\mathcal{F}_2|=g_1(k_1,k_2,n,t)g_2(k_2,n,t)$ if $\mathcal{F}_1$ and $\mathcal{F}_2$ are a pair of families given in ${\rm (ia)}$ or ${\rm (ib)},$ and  $|\mathcal{F}_1||\mathcal{F}_2|=g_3(k_1,k_2,n,t)$ if $\mathcal{F}_1$ and $\mathcal{F}_2$ are a pair of families given in ${\rm (iia)}$ or ${\rm (iib)}.$ Now suppose that $\mathcal{F}_1$ and $\mathcal{F}_2$ are maximal cross $t$-intersecting families, which are neither a pair of families given in ${\rm(ia)}$ and ${\rm(ib)}$, nor a pair of families given in ${\rm(iia)}$ and ${\rm(iib)}$. To prove the theorem, we only need to show that
\begin{align*}
|\mathcal{F}_1||\mathcal{F}_2|<\left\{\begin{array}{ll}
g_1(k_1,k_2,n,t)g_2(k_2,n,t),& \mbox{if}\ k_2\geq 2t+1,\\
g_3(k_1,k_2,n,t),& \mbox{if}\  k_2\leq 2t.
\end{array}\right.
\end{align*}
By (\ref{v-ineq-10}), it suffices to prove that
\begin{align}\label{F1F2upp}
|\mathcal{F}_1||\mathcal{F}_2|<g_1(k_1,k_2,n,t)g_2(k_2,n,t)\quad \mbox{or}\quad  |\mathcal{F}_1||\mathcal{F}_2|<g_3(k_1,k_2,n,t).
\end{align}
 We divide our proof into the following five cases.

\medskip

\noindent\textbf{Case 1.} $\tau_t(\mathcal{F}_1)=\tau_t(\mathcal{F}_2)=t$.

Assume that $X_1$ and $X_2$ are $t$-covers of $\mathcal{F}_1$ and $\mathcal{F}_2$ with dimensional $t$, respectively. By Lemma~\ref{v-FTt}, we have $X_1=X_2$, implying that $\mathcal{F}_1$ and $\mathcal{F}_2$ are trivial cross $t$-intersecting families, a contradiction.

\medskip

\noindent\textbf{Case 2.} $\tau_t(\mathcal{F}_1)=t$ and $\tau_t(\mathcal{F}_2)=t+1$.

Assume that $X$ is a $t$-cover of $\mathcal{F}_1$ with dimensional $t$, and $T$ is a $t$-cover of $\mathcal{F}_2$ with dimensional $t+1$.

Suppose that there exists a $(k_2+1)$-subspace $M$ of $V$ such that $F_2+X=M$ for each $F_2\in\mathcal{F}_2\setminus(\mathcal{F}_2)_X$. Then $\mathcal{F}_1=\mathcal{A}(k_1,t+1,X,M)$ and $\mathcal{F}_2=\mathcal{B}(k_2,X,M)$ by Lemma~\ref{lem3-1} (ii), and $|\mathcal{F}_1||\mathcal{F}_2|=g_1(k_1,k_2,n,t)g_2(k_2,n,t)$ by Lemma~\ref{lem-abcd}. Since $\mathcal{F}_1$ and $\mathcal{F}_2$ are not a pair of families given in ${\rm (ia)}$ or ${\rm (ib)},$ we have $k_2\leq 2t$. If $t+2\leq k_2\leq 2t$, by Lemmas~\ref{lem-upp-1} and \ref{v-inequ-3}, then
$|\mathcal{F}_1||\mathcal{F}_2|<g_3(k_1,k_2,n,t).$
If $k_2=t+1$, by Lemma~\ref{v-inequ-3-1}, then $|\mathcal{F}_1||\mathcal{F}_2|<g_3(k_1,k_2,n,t).$ Therefore, (\ref{F1F2upp}) holds.

Suppose that $T\subseteq F_1$ for each $F_1\in\mathcal{F}_1.$ Then $\mathcal{F}_1=\mathcal{C}(k_1,T)$, $\mathcal{F}_2=\mathcal{D}(k_2,t,T)$ and $|\mathcal{F}_1||\mathcal{F}_2|=g_3(k_1,k_2,n,t)$ by Lemmas~\ref{lem3-1} (iii) and \ref{lem-abcd} (ii). Since $\mathcal{F}_1$ and $\mathcal{F}_2$ are not a pair of families given in ${\rm (iia)}$ or ${\rm (iib)}$, we have $k_2\geq 2t+1$. It follows from Lemmas~\ref{lem-upp-1} and \ref{v-inequ-2} that
$|\mathcal{F}_1||\mathcal{F}_2|<g_1(k_1,k_2,n,t)g_2(k_2,n,t),$
and so (\ref{F1F2upp}) holds.

Suppose that there exist $F_{2,1},\ F_{2,2}\in\mathcal{F}_2\setminus(\mathcal{F}_2)_X$ such that $F_{2,1}+X\neq F_{2,2}+X$, and there exists $F_1\in\mathcal{F}_1$ such that $T\nsubseteq F_1$. It follows from Lemmas~\ref{lem3-1} (iv), \ref{lem-upp-1} and \ref{v-inequ-6}  that
$$
|\mathcal{F}_1||\mathcal{F}_2|\leq f_3(k_1,k_2,n,t)f_4(k_1,k_2,n,t)<g_1(k_1,k_2,n,t)g_2(k_2,n,t),
$$
and so (\ref{F1F2upp}) holds.

\medskip

\noindent\textbf{Case 3.} $\tau_t(\mathcal{F}_1)=t+1$ and $\tau_t(\mathcal{F}_2)=t$.

Assume that $T$ is a $t$-cover of $\mathcal{F}_1$ with dimensional $t+1$, and $X$ is a $t$-cover of $\mathcal{F}_2$ with dimensional $t$.

Suppose that there exists a $(k_1+1)$-subspace $M$ of $V$ such that $F_1+X=M$ for each $F_1\in\mathcal{F}_1\setminus(\mathcal{F}_1)_X$. Then $\mathcal{F}_1=\mathcal{B}(k_1,X,M)$ and $\mathcal{F}_2=\mathcal{A}(k_2,t+1,X,M)$ by Lemma~\ref{lem3-1} (ii), and $|\mathcal{F}_1||\mathcal{F}_2|=g_2(k_1,n,t)g_1(k_2,k_1,n,t)$ by Lemma~\ref{lem-abcd} (i). If $k_1>k_2$,  then $|\mathcal{F}_1||\mathcal{F}_2|<g_1(k_1,k_2,n,t)g_2(k_2,n,t)$ from Lemma~\ref{v-inequ-5}. Now assume that $k_1=k_2$. Since $\mathcal{F}_1$ and $\mathcal{F}_2$ are not a pair of families given in ${\rm (ib)},$ we have $k_2\leq 2t$, implying that $|\mathcal{F}_1||\mathcal{F}_2|=g_2(k_1,n,t)g_1(k_2,k_1,n,t)=g_1(k_1,k_2,n,t)g_2(k_2,n,t)<g_3(k_1,k_2,n,t)$ from (\ref{v-ineq-10}). Therefore (\ref{F1F2upp}) holds.

Suppose that $T\subseteq F_2$ for each $F_2\in\mathcal{F}_2.$ Then $\mathcal{F}_1=\mathcal{D}(k_1,t,T)$, $\mathcal{F}_2=\mathcal{C}(k_2,T)$ and $|\mathcal{F}_1||\mathcal{F}_2|=g_3(k_2,k_1,n,t)$ by Lemmas~\ref{lem3-1} (iii) and \ref{lem-abcd} (ii). If $k_1>k_2$, then $|\mathcal{F}_1||\mathcal{F}_2|<g_3(k_1,k_2,n,t)$ from Lemma~\ref{v-inequ-1}. Now assume that $k_1=k_2$. Since $\mathcal{F}_1$ and $\mathcal{F}_2$ are not a pair of families given in ${\rm (iib)},$ we have $k_2\geq 2t+1$, implying that $|\mathcal{F}_1||\mathcal{F}_2|=g_3(k_2,k_1,n,t)=g_3(k_1,k_2,n,t)<g_1(k_1,k_2,n,t)g_2(k_2,n,t)$ by (\ref{v-ineq-10}). Therefore (\ref{F1F2upp}) holds.

Suppose that  there exist $F_{1,1},\ F_{1,2}\in\mathcal{F}_1\setminus(\mathcal{F}_1)_X$ such that $F_{1,1}+X\neq F_{1,2}+X$, and there exists $F_2\in\mathcal{F}_2$ such that $T\nsubseteq F_2$. It follows from Lemmas~\ref{lem3-1} (iv), \ref{lem-upp-1} and \ref{v-inequ-6} that
$$
|\mathcal{F}_1||\mathcal{F}_2|\leq f_3(k_2,k_1,n,t)f_4(k_2,k_1,n,t)<g_1(k_1,k_2,n,t)g_2(k_2,n,t),
$$
and hence (\ref{F1F2upp}) holds.

\medskip

\noindent\textbf{Case 4.} $\tau_t(\mathcal{F}_1)=t$ and $\tau_t(\mathcal{F}_2)\geq t+2$, or $\tau_t(\mathcal{F}_1)\geq t+2$ and $\tau_t(\mathcal{F}_2)=t$.

Assume that $\tau_t(\mathcal{F}_1)=t$ and $\tau_t(\mathcal{F}_2)\geq t+2$. By Corollary~\ref{v-cor} (i), Lemmas~\ref{lem-upp-1} and \ref{v-inequ-8}, we have that
$$
|\mathcal{F}_1||\mathcal{F}_2|\leq {n-t\brack k_2-t}f_2(t+2,k_2,k_1,n,t)<g_1(k_1,k_2,n,t)g_2(k_2,n,t),
$$
and (\ref{F1F2upp}) holds.

Assume that  $\tau_t(\mathcal{F}_1)\geq t+2$ and $\tau_t(\mathcal{F}_2)=t$. By Corollary~\ref{v-cor}  (i), we have that
$$
|\mathcal{F}_1||\mathcal{F}_2|\leq {n-t\brack k_1-t}f_2(t+2,k_1,k_2,n,t),
$$
implying that (\ref{F1F2upp}) holds from Lemmas~\ref{lem-upp-1} and \ref{v-inequ-9}.

\medskip

\noindent\textbf{Case 5.} $\tau_t(\mathcal{F}_1)\geq t+1$ and $\tau_t(\mathcal{F}_2)\geq t+1$.			 

Assume that $\tau_t(\mathcal{F}_1)=m_1$ and $\tau_t(\mathcal{F}_2)=m_2$. From Corollary~\ref{v-cor} (ii), Lemmas \ref{lem-upp-1} and \ref{v-inequ-4}, we obtain
$$
|\mathcal{F}_1||\mathcal{F}_2|<f_1(k_1,k_2,n,t)f_1(k_2,k_1,n,t)<g_1(k_1,k_2,n,t)g_2(k_2,n,t),
$$
and (\ref{F1F2upp}) holds.    $\qed$

\section{Non-trivial $r$-wise $t$-intersecting families for vector spaces}

In this section, we determine the structure of the non-trivial maximal $r$-wise $t$-intersecting families for subspaces of $V$, and give a stability result about these families.
\subsection{Some properties}
\begin{lemma}\label{v-lem-1}
	Let $\mathcal{F}\subseteq {V\brack k}$ be a non-trivial $r$-wise $t$-intersecting family. Suppose that $S$ is a subspace of $V$ such that $\dim(S\cap F_1\cap F_2\cap\cdots \cap F_{r-1})\geq t$ for all $F_1,F_2,\ldots,F_{r-1}\in\mathcal{F}$. Then the following hold.
	\begin{enumerate}[\rm(i)]
		\item We have $\dim(S)\geq t+r-1.$
		\item Suppose that $A_1,A_2,\ldots,A_m$ are elements in $\mathcal{F}$ with $m\leq r-1$. Then $\dim(S\cap A_1\cap A_2\cap \cdots\cap A_m)\geq t+r-m-1.$
	\end{enumerate}
\end{lemma}
\proof (i)\quad Since $\mathcal{F}$ is non-trivial, for each subspace $T$ of $V$ with $\dim(T)\geq t$, there exists an element $A\in\mathcal{F}$ such that $\dim(T\cap A)\leq \dim(T)-1$. It is clear that $\dim(S)\geq t$.  Then, there exists $A_1\in\mathcal{F}$ such that $t\leq \dim(S\cap A_1)\leq \dim(S)-1$, and (i) holds for $r=2$. Now suppose that $r\geq 3$, by induction, there exist $A_2,\ldots, A_{r-1}\in\mathcal{F}$ such that $t\leq\dim(S\cap A_1\cap\cdots\cap A_i)\leq \dim(S\cap A_1\cap\cdots\cap A_{i-1})-1$ for each $i\in\{2,3,\ldots,r-1\}$ due to $\dim(S\cap F_1\cap F_2\cap\cdots \cap F_{r-1})\geq t$ for all $F_1,F_2,\ldots,F_{r-1}\in\mathcal{F}$. Therefore
$$
\dim(S)\geq\dim(S\cap A_1)+1\geq \cdots\geq \dim(S\cap A_1\cap \cdots\cap A_{r-1})+r-1\geq t+r-1
$$
 as required.

(ii)\quad Observe that it is clear if $m=r-1$. Suppose that $m<r-1$ in the following. Since $\dim((S\cap A_1\cap A_2\cap\cdots\cap A_m)\cap F_1\cap\cdots\cap F_{r-1-m})\geq t$ for all $F_1,F_2,\ldots,F_{r-1-m}\in\mathcal{F}$, and $\mathcal{F}$ is a non-trivial $(r-m)$-wise $t$-intersecting family, by (i), we have $\dim(S\cap A_1\cap\cdots\cap A_m)\geq t+r-m-1$ as required.    $\qed$
\begin{cor}
If $r>k-t+1,$ then there does not exist a non-trivial $r$-wise $t$-intersecting  family of $k$-subspaces of $V$.
\end{cor}
\proof Suppose that there exists a non-trivial $r$-wise $t$-intersecting family $\mathcal{F}$. Then $\dim(F_1\cap \cdots\cap F_r)\geq t$ for all $F_1, F_2, \ldots,F_r\in\mathcal{F}$. By Lemma~\ref{v-lem-1} (i), we have $k=\dim(F_1)\geq t+r-1$, a contradiction. Hence, the desired result holds.  $\qed$
\begin{cor}\label{v-lem-2}
	Suppose that $B_1, B_2,\ldots,B_d$ are elements in the non-trivial $r$-wise $t$-intersecting family $\mathcal{F}\subseteq{V\brack k}$ with $d\leq r$.
	Then $\dim(B_1\cap B_2\cap\cdots\cap B_d)\geq t+r-d$. Moreover, $\mathcal{F}$ is a non-trivial $(t+r-2)$-intersecting family.
\end{cor}
\proof Set $S=B_1$ in Lemma~\ref{v-lem-1} (ii), and then the former part of this lemma holds. In particular, set $d=2$, and then we have $\mathcal{F}$ is $(t+r-2)$-intersecting.   $\qed$

\subsection{The proof of Theorem~\ref{v-main-3}}

Let $d$, $k$ and $n$ be positive integers with $n\geq 2k$. Let $X,M$ and $C$ be subspaces of $V$ such that $X\subseteq M\subseteq C,$ $\dim(X)=d$, $\dim(M)=k$ and $\dim(C)=c$, where $c\in\{k+1,k+2,\ldots, 2k-d,n\}$. Define
\begin{align}
	\mathcal{H}_2(X,M,C) &= \mathcal{E}_1(X,M) \cup \mathcal{E}_2(X, M,C) \cup \mathcal{E}_3(X, M, C)\label{eq:H2},
\end{align}
where
\begin{eqnarray*}
	\mathcal{E}_1(X,M) &=&\left\{F\in{V\brack k}\mid X\subseteq F,\ \dim(F\cap M)\geq d+1\right\},\\
	\mathcal{E}_2(X, M,C) &=&\left\{F\in{V\brack k}\mid F\cap M=X,\ \dim(F\cap C)=c-k+d\right\},\\
	\mathcal{E}_3(X, M, C) &=&\left\{F\in{C\brack k}\mid \dim(F\cap X)=d-1,\ \dim(F\cap M)=k-1\right\}. 
\end{eqnarray*}

From \cite[Remark~1 and Lemma~2.5]{Cao-vec}, the following hold.
\begin{obs}\label{v-obs1} Let $h_1(d,k,n)$ and $h_2(d,k,n)$ be as in (\ref{eqq1}) and (\ref{eqq2}).
\begin{enumerate}[{\rm(i)}]
\item	$\mathcal{H}_2(X,M,C) = \mathcal{A}(k,d+1,X,C)\cup{C\brack k}$ and $|\mathcal{H}_2(X,M,C)|=h_1(d,k,n)$ if $C$ satisfies $\dim(C) = k+1$.
\item $\mathcal{H}_2(X,M,V)=\mathcal{E}_1(X,M)\cup \mathcal{E}_3(X, M, V)$ and $|\mathcal{H}_2(X,M,V)|=h_2(d,k,n)$.
\item $\mathcal{H}_2(X,M,V)=\mathcal{D}(k,d+1,M)$ if $d$ and $k$ satisfy $d=k-2$.
\end{enumerate}
\end{obs}

 The following two theorems given in \cite{Cao-vec} are essential to prove Theorem~\ref{v-main-3}.

\begin{thm}{\rm (\cite[Theorem~1.1]{Cao-vec})}\label{v-non-trivial-0}
	Let $n, k$ and $t$ be positive integers with $t\leq k-2$ and $2k+t+\min\{4, 2t\}\leq n$. If $\mathcal{F}\subseteq {V\brack k}$ is a maximal non-trivial $t$-intersecting family and
	$$|\mathcal{F}| \geq {k-t\brack 1}{n-t-1\brack k-t-1}-q{k-t\brack 2}{n-t-2\brack k-t-2},$$ then one of the following holds:
	\begin{itemize}
		\item[{\rm(i)}] $\mathcal{F}=\mathcal{H}_2(X,M,C)$ for some $t$-subspace $X$, $k$-subspace $M$ and $c$-subspace $C$ of $V$ with $X\subseteq M\subseteq C$ and $c\in\{k+1,k+2,\ldots, 2k-t,n\}$;
		\item[{\rm(ii)}] $\mathcal{F}=\mathcal{D}(k,t+1,Z)$ for some $(t+2)$-subspace $Z$ of $V$, and $\frac{k}{2}-1\leq t\leq k-2$.
	\end{itemize}
\end{thm}

\begin{thm}{\rm (\cite[Theorem~1.2]{Cao-vec})}\label{v-non-trivial-1}
	Let $n, k$ and $t$ be positive integers with $t\leq k-2$ and $2k+t+\min\{4, 2t\}\leq n$, and let $\mathcal{F}\subseteq {V\brack k}$ be a non-trivial $t$-intersecting family. Then the following hold.
	\begin{itemize}
		\item[{\rm(i)}]  If $1\leq t\leq \frac{k}{2}-1$, then
		$$
		|\mathcal{F}|\leq {n-t\brack k-t}-q^{(k+1-t)(k-t)}{n-k-1\brack k-t}+q^{k+1-t}{t\brack 1},
		$$
		and equality holds if and only if $\mathcal{F}=\mathcal{A}(k,t+1,X,M)\cup{M\brack k}$ for some $t$-subspace $X$ and $(k+1)$-subspace $M$ of $V$ with $X\subset M$.
		\item[{\rm(ii)}] If $\frac{k}{2}-1<t\leq k-2,$ then
		$$
		|\mathcal{F}|\leq {t+2\brack 1}{n-t-1\brack k-t-1}-q{t+1\brack 1}{n-t-2\brack k-t-2},
		$$
		and equality holds if and only if $\mathcal{F}=\mathcal{D}(k,t+1,Z)$ for some $(t+2)$-subspace $Z$ of $V$, or $(t, k) = (1, 3)$ and $\mathcal{F}=\mathcal{A}(k,t+1,X,M)\cup{M\brack k}$ for some $1$-subspace $X$ and $4$-subspace $M$ of $V$ with $X\subset M$.
	\end{itemize}
\end{thm}
\begin{lemma} \label{v-family}
Let $X\in{V\brack t+r-2}$ and $M\in{V\brack k+1}$ with $X\subset M$, and  $Z\in{V\brack t+r}$. Then both $\mathcal{A}(k,t+r-1,X,M)\cup{M\brack k}$ and $\mathcal{D}(k,t+r-1,Z)$  are non-trivial $r$-wise $t$-intersecting families.
\end{lemma}
\proof Firstly, we claim that for positive integers $\ell$ and $m$ with $m\leq \ell,$ if $E_1,\ldots,E_m$ are $\ell$-subspaces of a space $E$ with $\dim(E)=\ell+1,$ then $\dim(E_1\cap E_2\cap\cdots\cap E_m)\geq\ell-m+1$. Indeed,
\begin{align*}
	\dim(E_1\cap E_2\cap\cdots\cap E_m)\geq&\dim(E_1\cap E_2\cap\cdots\cap E_{m-1})+\dim(E_m)-\dim(E)\\
	=&\dim(E_1\cap E_2\cap\cdots\cap E_{m-1})-1\\
	\geq&\dim(E_1\cap E_2\cap\cdots\cap E_{m-2})-2\\
	&\vdots\\
	\geq& \dim(E_1)-(m-1)\\
	=&\ell-m+1.
\end{align*}

Let $F_1,F_2,\ldots,F_r\in\mathcal{A}(k,t+r-1,X,M)\cup{M\brack k}.$ If $F_i\in{M\brack k}$ for each $i\in\{1,2,\ldots,r\}$, then $\dim(F_1\cap\cdots\cap F_r)\geq t$ due to the claim above and $r\leq k-t+1$. Otherwise, without loss of generality, suppose that there exists $h\in\{1,2,\ldots,r\}$ such that
$$
F_i\in\left\{\begin{array}{ll} {M\brack k}, & i\leq h,\\\mathcal{A}(k,t+r-1,X,M)\setminus{M\brack k},& \mbox{otherwise.}\end{array}\right.
$$
Then, if $h=r-1$, we have
\begin{align*}
	\dim(F_1\cap\cdots\cap F_r)\geq& \dim((F_1\cap\cdots\cap F_{r-1})\cap (F_r\cap M))\\
	\geq&\dim(F_1\cap\cdots\cap F_{r-1})+\dim (F_r\cap M)-\dim(M)\\
	\geq&(k-r+2)+(t+r-1)-(k+1)\\
	\geq&t;
\end{align*}
if $h\leq r-2$, we have
\begin{align*}
	\dim(F_1\cap\cdots\cap F_r)\geq& \dim((F_1\cap\cdots\cap F_h)\cap (F_{h+1}\cap\cdots\cap F_r\cap X))\\
	\geq&\dim(F_1\cap\cdots\cap F_h)+\dim(F_{h+1}\cap\cdots\cap F_r\cap X)-\dim(M)\\
	\geq&(k-h+1)+(t+r-2)-(k+1)\\
	\geq&t.
\end{align*}

Let $F_1,F_2,\ldots,F_r\in\mathcal{D}(k,t+r-1,Z).$ If $Z\subseteq F_i$ for each $i\in\{1,2,\ldots,r\}$, then $\dim(F_1\cap\cdots\cap F_r)\geq t$. Otherwise, without loss of generality, suppose that there exists $h\in\{1,2,\ldots,r\}$ such that
$$
\dim(F_i\cap Z)=\left\{\begin{array}{ll} t+r-1, & i\leq h,\\t+r,& \mbox{otherwise.}\end{array}\right.
$$
Then, by the claim above, we have
\begin{align*}
	\dim(F_1\cap\cdots\cap F_r)\geq&\dim((F_1\cap Z)\cap\cdots\cap (F_r\cap Z))\\
	=&\dim((F_1\cap Z)\cap\cdots\cap (F_h\cap Z))\\
	\geq&t+r-1-h+1\\
	\geq&t.
\end{align*}

Hence, both $\mathcal{A}(k,t+r-1,X,M)\cup{M\brack k}$ and $\mathcal{D}(k,t+r-1,Z)$ are $r$-wise $t$-intersecting. It is clear that these two families are non-trivial, and therefore the lemma holds.  $\qed$
\begin{lemma}
		Suppose that $\mathcal{F}\subseteq{V\brack k}$ is a maximal non-trivial $r$-wise $t$-intersecting family with $t+r=k+1$ and $k+1\leq n$. Then $\mathcal{F}={M\brack k}$ for some $(k+1)$-subspace $M$.
\end{lemma}
\proof It is known that $\mathcal{F}$ is $(k-1)$-intersecting by Corollary~\ref{v-lem-2}. By \cite[Remark (ii) in Section 9.3]{ABE}, we have $\mathcal{F}\subseteq{M\brack k}$ for some $(k+1)$-subspace $M$. Observe that ${M\brack k}$ is $r$-wise $t$-intersecting due to $t+r=k+1$ and the claim in the proof of Lemma~\ref{v-family}. Then we obtain $\mathcal{F}={M\brack k}$ due to the maximality of $\mathcal{F}$.   $\qed$
\begin{lemma}\label{v-wise-sub}
Let $n$, $k$, $t$ and $r$ be positive integers with $r\geq 3$, $t+r-2\leq k-2$ and $2k+t+r+2\leq n$. Let $\mathcal{F}$ be a non-trivial $r$-wise $t$-intersecting subfamily of $\mathcal{H}_2(X,M,C)$, where $X\subseteq M\subseteq C,$ $\dim(X)=t+r-2$, $\dim(M)=k$ and $\dim(C)=c\in\{k+2,k+3,\ldots, 2k-t-r+2\}$. If $|\mathcal{F}|> h_2(t+r-2,k,n),$ then $\mathcal{F}\subseteq\mathcal\mathcal{A}(k,t+r-1,X,M^\prime)\cup{M^\prime\brack k}$ for some $(k+1)$-subspace $M^{\prime}$ of $C$ with $M\subseteq M^{\prime}$.
\end{lemma}
\proof Suppose for the contrary that $\mathcal{F}\nsubseteq \mathcal{A}(k,t+r-1,X,M^\prime)\cup{M^\prime\brack k}$ for each $(k+1)$-subspace $M^\prime$ of $C$ with $M\subseteq M^\prime.$ Let $\mathcal{E}_1(X,M),$ $\mathcal{E}_2(X, M,C)$ and $\mathcal{E}_3(X, M, C)$ be as in (\ref{eq:H2}). Observe that all the subspaces in $\mathcal{E}_1(X,M)$ and $\mathcal{E}_2(X, M,C)$ are containing $X$. Since $\mathcal{F}$ is non-trivial, we have $\mathcal{F}\cap \mathcal{E}_3(X, M, C)\neq\emptyset.$

\textbf{Claim~1.} For each subspace $S$ of $X$ with $\dim(S)\geq t$, there exists $F\in \mathcal{F}\cap \mathcal{E}_3(X, M, C)$ such that $\dim(S\cap F\cap X)=\dim(S)-1.$

Indeed, since $\mathcal{F}$ is non-trivial and $\dim(S)\geq t$, there exists $F\in\mathcal{F}\cap \mathcal{E}_3(X, M, C)$ such that $S\nsubseteq F.$ From $\dim(F\cap X)=\dim(X)-1$, observe that $S+(F\cap X)=X$ and
$$
\dim(S\cap F\cap X)=\dim(S)+\dim(F\cap X)-\dim(S+(F\cap X))=\dim(S)-1.
$$
Hence, Claim 1 holds.

\medskip

\textbf{Claim~2.} There exist $F_1,\;F_2\in\mathcal{F}\cap \mathcal{E}_3(X, M, C)$ satisfying $F_1\cap X\neq F_2\cap X$ and $F_2\nsubseteq F_1+M.$

Choose $G_1\in \mathcal{F}\cap \mathcal{E}_3(X, M, C)$. We have $G_1+M\subseteq C,$ and $\dim(G_1+M)=k+1$ from $\dim(G_1\cap M)=k-1.$ By the construction of $\mathcal{E}_2(X,M,C)$, observe that $\mathcal{E}_2(X,M,C)\subseteq\mathcal{E}_2(X,M,G_1+M)$, implying that
$$
\mathcal{E}_1(X,M)\cup \mathcal{E}_2(X,M,C)\subseteq \mathcal{E}_1(X,M)\cup \mathcal{E}_2(X,M,G_1+M)\subseteq\mathcal{A}(k,t+r-1, X, G_1+M).
$$
Then $\mathcal{F}\cap(\mathcal{E}_1(X,M)\cup \mathcal{E}_2(X,M,C))\subseteq\mathcal{A}(k,t+r-1, X, G_1+M).$ Since $\mathcal{F}\nsubseteq \mathcal{A}(k,t+r-1,X,G_1+M)\cup{G_1+M\brack k},$ we have $\mathcal{F}\cap\mathcal{E}_3(X,M,C)\nsubseteq {G_1+M\brack k}$, implying that
there exists $G_2\in\mathcal{F}\cap\mathcal{E}_3(X,M,C)$ such that $G_2\nsubseteq G_1+M.$

If $G_2\cap X\neq G_1\cap X$, setting $F_1=G_1$ and $F_2=G_2$, then $F_1$ and $F_2$ are the required subspace.

Now assume that $G_3\cap X=G_1\cap X$ for each $G_3\in\mathcal{F}\cap\mathcal{E}_3(X,M,C)$ with $G_3\nsubseteq G_1+M$. Observe that $X\subseteq F$ for each $F\in\mathcal{F}\cap(\mathcal{E}_1(X,M)\cup\mathcal{E}_2(X,M,C))$, and $\dim(F\cap X)=t+r-3\geq t$ for each $F\in\mathcal{F}\cap\mathcal{E}_3(X,M,C)$. Since $\mathcal{F}$ is non-trivial, there exists $G_4\in\mathcal{F}\cap\mathcal{E}_3(X,M,C)$ such that $G_4\cap X\neq G_1\cap X=G_2\cap X$. It follows from the assumption that $G_4\in{G_1+M\brack k}$, and so $G_4+M=G_1+M$ due to $\dim(G_4\cap M)=k-1$ and $\dim(G_1+M)=k+1$. Set $F_1=G_4$ and $F_2=G_2$. We have that $F_1$ and $F_2$ are the required subspace, and Claim 2 holds.

\medskip

By Claim~2, let $F_1$ and $F_2$ be the subspace in $\mathcal{F}\cap\mathcal{E}_3(X,M,C)$ with  $F_1\cap X\neq F_2\cap X$ and $F_2\nsubseteq F_1+M.$ Then $\dim(F_1\cap F_2\cap X)=\dim(F_1\cap X)+\dim(F_2\cap X)-\dim((F_1\cap M)+(F_2\cap X))=t+r-4$. Using Claim $1$ repeatedly, we can get $F_1,\;F_2,\ldots,F_{r-1}\in\mathcal{F}\cap\mathcal{E}_3(X,M,C)$ satisfying $F_1\cap X\neq F_2\cap X,$ $F_2\nsubseteq F_1+M$  and $\dim(X\cap F_1\cap F_2\cap\cdots\cap F_{r-1})=t-1.$

Since $\mathcal{F}\subseteq \mathcal{H}_2(X,M,C)$ and $|\mathcal{F}|> h_2(t+r-2,k,n)$, by Observation~\ref{v-obs1} (ii),
we have
\begin{align*}
	|\mathcal{F}\cap\mathcal{E}_2(X, M,C)|\geq& |\mathcal{F}|-|\mathcal{E}_1(X,M)|-|\mathcal{E}_3(X, M, C)|\\
	>& |\mathcal{H}_2(X,M,V)|-|\mathcal{E}_1(X,M)|-|\mathcal{E}_3(X, M, C)|\\
	=&|\mathcal{E}_3(X, M, V)|-|\mathcal{E}_3(X, M, C)|>0,
\end{align*}
which implies that $\mathcal{F}\cap\mathcal{E}_2(X, M,C)\neq\emptyset$.

Let $F\in\mathcal{F}\cap\mathcal{E}_2(X, M,C).$  Then $\dim(F\cap F_1\cap F_2\cap\cdots\cap F_{r-1})\geq t$ because $\mathcal{F}$ is $r$-wise $t$-intersecting. Since $\dim(X\cap F_1\cap F_2\cap\cdots\cap F_{r-1})=t-1$, we have $F\cap F_1\cap F_2\neq X\cap F_1\cap F_2$. It is clear that $F\cap F_1\cap F_2\supseteq X\cap F_1\cap F_2.$ Then there exists $y\in(F\cap F_1\cap F_2)\setminus (X\cap F_1\cap F_2)$. By the construction of $\mathcal{E}_2(X, M,C)$, note that $F\cap M=X$. It follows that $y\notin M\cap F_1$ and $y\notin M\cap F_2$. Hence $F_1=(M\cap F_1)+\langle y\rangle$ and $F_2=(M\cap F_2)+\langle y\rangle$, implying that $F_2\subseteq M+\langle y\rangle=F_1+M,$ a contradiction.  $\qed$

\noindent{\em \textbf{Proof of Theorem~\ref{v-main-3}}}\quad   Since $\mathcal{F}\subseteq{V\brack k}$ is a non-trivial $r$-wise $t$-intersecting family, by Corollary~\ref{v-lem-2}, $\mathcal{F}$ is a non-trivial $(t+r-2)$-intersecting family. It follows from $r\geq 3$, $t+r-2\leq k-2$ and $2k+t+r+3\leq n$ that $2k+(t+r-2)+\min\{4,2t+2r-4\}\leq n$ and $(t+r-2,k)\neq(1,3).$ Then by Theorem~\ref{v-non-trivial-1} and Lemma~\ref{v-family}, we have that each former part of (i) and (ii) in Theorem~\ref{v-main-3} holds.

Suppose that $|\mathcal{F}|>h_2(t+r-2,k,n)$ if $t+r-2<k-2$, and $|\mathcal{F}|>h_1(t+r-2,k,n)$ if $t+r-2=k-2$. From \cite[Lemmas 2.6 (i) and 2.7 (i)]{Cao-vec}, we have that
$$
|\mathcal{F}|>{k-t-r+2\brack 1}{n-t-r+1\brack k-t-r+1}-q{k-t-r+2\brack 2}{n-t-r\brack k-t-r}.
$$
Then by Theorem~\ref{v-non-trivial-0}, one of the following holds:
\begin{enumerate}[{\rm (i)}]
\item $\mathcal{F}\subseteq\mathcal{H}_2(X,M,C)$ for some $X\in{V\brack t+r-2}$, $M\in{V\brack k}$ and $C\in{V\brack c}$ with $X\subseteq M\subseteq C$ and $c\in\{k+1,k+2,\ldots, 2k-t-r+2,n\}$;
\item $\mathcal{F}\subseteq\mathcal{D}(k,t+r-1,Z)$ for some  $Z\in{V\brack t+r}$, and $\frac{k}{2}-1\leq t+r-2\leq k-2$.	
\end{enumerate}

For $t+r-2=k-2$, recall that $h_2(t+r-2,k,n)> h_1(t+r-2,k,n)$  due to \cite[Lemma~2.6(iii)]{Cao-vec}, and $\mathcal{H}_2(X,M,V)=\mathcal{D}(k,k-1,M)$ by Observation~\ref{v-obs1} (iii). Hence, by Lemmas~\ref{v-family} and \ref{v-wise-sub}, the latter parts of (i) and (ii) in Theorem~\ref{v-main-3} hold.   $\qed$

\section{Some inequalities}
In this section, we prove some inequalities used in this paper. In the following lemmas, we always assume that $n$, $k_1$, $k_2$ and $t$ are positive integers with $n\geq k_1+k_2+t+3$ and $k_1\geq k_2\geq t+1$, and let $f_j(k,\ell,n,t)$ and $g_i(k,\ell,n,t)$ be as in (\ref{eqq3})--(\ref{eqq4}), (\ref{eqq5})--(\ref{eqq11}). For convenience, we set $\prod_{i=a}^{a-1}b_i=1$ in the following. Write
{\small \begin{align}
	g_4(k,\ell,n,t)=&{\ell-t+1\brack 1}{n-t-1\brack k-t-1}-q{\ell-t+1\brack 2}{n-t-2\brack k-t-2},\label{eqq7}\\
	g_5(k,\ell,n,t)=& {\ell-t+1\brack 1}{n-t-1\brack k-t-1}-q^{(\ell-t-1)(k-t-2)+1}{n-\ell-1\brack k-t-2}{\ell+1-t\brack 2}.\label{eqq9}
\end{align}}
\begin{lemma}\label{lem-upp-1}
	If $\min\{k,\ell\}\geq t+1,$ then $g_4(k,\ell,n,t)\leq g_1(k,\ell,n,t)\leq g_5(k,\ell,n,t)$.
\end{lemma}
\proof If $k=t+1,$ then $g_1(k,\ell,n,t)={n-t\brack 1}-q^{(\ell+1-t)}{n-\ell-1\brack 1}={\ell+1-t\brack 1}.$ It is routine to check that the result holds for this case. Now suppose that $k\geq t+2.$ Let $M$ be a $(\ell+1)$-subspace of $V$, and $X$ a $t$-subspace of $M$. For each $i\in\{t,t+1,\ldots,\ell+1\}$, set
{\small\begin{align*}
	\mathcal{A}_i(X,M)=&\left\{F\in{V\brack k}\mid X\subseteq F,\ \dim(F\cap M)=i\right\},\\
	\mathcal{L}_i(X,M)=&\left\{(I,F)\in{V\brack i}\times{V\brack k}\mid X\subseteq I\subseteq M,\ I\subseteq F\right\}.
\end{align*}}
Observe that $$\sum_{j=t+1}^{\min\{k,\ell+1\}}|\mathcal{A}_j(X,M)|=|\mathcal{A}(k,t+1,X,M)|=g_1(k,\ell,n,t)$$ from Lemma~\ref{lem-abcd} (i). By double counting $|\mathcal{L}_i(X,M)|$, we have that
{\small\begin{align*}
	|\mathcal{L}_i(X,M)|=\sum_{j=i}^{\min\{k,\ell+1\}}|\mathcal{A}_j(X,M)|{j-t\brack i-t}={\ell+1-t\brack i-t}{n-i\brack k-i}.
\end{align*}}
Then we obtain
{\small\begin{align*}
	|\mathcal{L}_{t+1}(X,M)|-q|\mathcal{L}_{t+2}(X,M)|=&{\ell-t+1\brack 1}{n-t-1\brack k-t-1}-q{\ell-t+1\brack 2}{n-t-2\brack k-t-2}\\
	 =&\sum_{j=t+1}^{\min\{k,\ell+1\}}|\mathcal{A}_j(X,M)|+\sum_{j=t+2}^{\min\{k,\ell+1\}}|\mathcal{A}_j(X,M)|\left({j-t\brack 1}-1-q{j-t\brack 2}\right)\\
	\leq&\sum_{j=t+1}^{\min\{k,\ell+1\}}|\mathcal{A}_j(X,M)|=g_1(k,\ell,n,t)
\end{align*}}
from ${j-t\brack 1}-1-q{j-t\brack 2}=q{j-t-1\brack 1}-q{j-t\brack 2}\leq 0$, and $g_1(k,\ell,n,t)\geq g_4(k,\ell,n,t)$ holds.

Since
{\small\begin{align*}
|\mathcal{L}_{t+1}(X,M)|=\sum_{j=t+1}^{\min\{k,\ell+1\}}|\mathcal{A}_j(X,M)|+\sum_{j=t+2}^{\min\{k,\ell+1\}}|\mathcal{A}_j(X,M)|\cdot\left({j-t\brack 1}-1\right),
\end{align*}}
and $|\mathcal{A}_{t+2}(X,M)|=N^\prime(t,t;k,t+2;n,n-\ell-1)$, we have
{\small\begin{align*}
	 g_1(k,\ell,n,t)=&|\mathcal{L}_{t+1}(X,M)|-\sum_{j=t+2}^{\min\{k,\ell+1\}}|\mathcal{A}_j(X,M)|\cdot\left({j-t\brack 1}-1\right)\\
\leq&|\mathcal{L}_{t+1}(X,M)|-|\mathcal{A}_{t+2}(X,M)|\cdot\left({2\brack 1}-1\right)\\
	=& {\ell-t+1\brack 1}{n-t-1\brack k-t-1}-q^{(\ell-t-1)(k-t-2)+1}{n-\ell-1\brack k-t-2}{\ell+1-t\brack 2},
\end{align*}}
and $g_1(k,\ell,n,t)\leq g_5(k,\ell,n,t)$ holds.     $\qed$

\begin{lemma}\label{v-inequ-2}
If $k_2\geq 2t+1,$ then $g_4(k_1,k_2,n,t)g_2(k_2,n,t)>g_3(k_1,k_2,n,t)$.	
\end{lemma}
\proof
By Lemma~\ref{lem1-1-1} (i) and (ii), observe that
{\small\begin{align*}
g_2(k_2,n,t)={n-t\brack k_2-t}+q^{k_2+1-t}{t\brack 1}>{n-t\brack k_2-t},
\end{align*}
{\small\begin{align*}
g_4(k_1,k_2,n,t)=&{k_2-t+1\brack 1}{n-t-1\brack k_1-t-1}-q{k_2-t+1\brack 2}{n-t-2\brack k_1-t-2}\\
=& {k_2-t\brack 1}{n-t-1\brack k_1-t-1}\left( \frac{q^{k_2-t+1}-1}{q^{k_2-t}-1}-\frac{q(q^{k_2-t+1}-1)}{q^2-1}\cdot\frac{q^{k_1-t-1}-1}{q^{n-t-1}-1}\right)\\
>& {k_2-t\brack 1}{n-t-1\brack k_1-t-1}\left( q-q^{1+(k_2-t)+k_1-n}\right),
\end{align*}
 and
{\small\begin{align*}
g_3(k_1,k_2,n,t)=&{n-t-1\brack k_1-t-1}\left(q^{k_2-t}{n-t-1\brack k_2-t}{t+1\brack 1}+{n-t-1\brack k_2-t-1}\right)\\
<&{n-t-1\brack k_1-t-1}{n-t\brack k_2-t}{t+1\brack 1}.
\end{align*}}
Since $n\geq k_1+k_2+t+3$ and $k_2\geq 2t+1$, we have $q-q^{1+(k_2-t)+k_1-n}>1$ and ${k_2-t\brack 1}\geq {t+1\brack 1}$. Then the desired result follows. $\qed$

\begin{lemma}\label{v-inequ-3}
If $t+2\leq k_2\leq 2t,$ then $g_5(k_1,k_2,n,t)g_2(k_2,n,t)<g_3(k_1,k_2,n,t)$.	
\end{lemma}
\proof Since $n\geq k_1+k_2+t+3$ and $k_1\geq k_2\geq t+2$, we have $(k_2-t-1)(n-k_2)\geq (k_2-t-1)(k_1+t+1)+ (k_2-t-1)\cdot2\geq 2t+3+ 2(k_2-t-1)\geq 2k_2+1\geq  2k_2-t+2$, and
{\small\begin{align*}
 &	{n-t-1\brack k_1-t-1}{n-t-1\brack k_2-t-1}-	{k_2-t+1\brack 1}{n-t-1\brack k_1-t-1}\cdot q^{k_2+1-t}{t\brack 1}\\
> &	{n-t-1\brack k_1-t-1}\left(q^{(k_2-t-1)(n-k_2)}-q^{2(k_2-t+1)+t}\right)\geq 0
\end{align*}}
by Lemma~\ref{lem1-1-1}(iii). Thus
{\small\begin{align}\label{v-ineq-11}
	&g_3(k_1,k_2,n,t)-g_5(k_1,k_2,n,t)g_2(k_1,k_2,n,t)\nonumber\\
	>& {n-t-1\brack k_1-t-1}\cdot q^{k_2-t}{t+1\brack 1}{n-t-1\brack k_2-t}
-{k_2-t+1\brack 1}{n-t-1\brack k_1-t-1}{n-t\brack k_2-t}\nonumber\\
&+q^{(k_2-t-1)(k_1-t-2)+1}{n-k_2-1\brack k_1-t-2}{k_2+1-t\brack 2}{n-t\brack k_2-t}.
\end{align}}

When $t+2\leq k_2\leq 2t-1$, we have
{\small\begin{align*}
	&\frac{{n-t-1\brack k_1-t-1}\cdot q^{k_2-t}{t+1\brack 1}{n-t-1\brack k_2-t}}{{k_2-t+1\brack 1}{n-t-1\brack k_1-t-1}{n-t\brack k_2-t}}=\frac{(q^{t+1}-1)(q^{n-t}-q^{k_2-t})}{(q^{k_2-t+1}-1)(q^{n-t}-1)}\geq\frac{q^{t+1}-1}{q^{t}-1}\cdot\frac{q^{n-t}-q^{t-1}}{q^{n-t}-1}\\
	&~~~~~~~~~~ >\left(q-\frac{q-1}{q^t-1}\right)\left(1-\frac{q^{t-1}-1}{q^{n-t}-1}\right)>q-\frac{q-1}{q^t-1}-\frac{q(q^{t-1}-1)}{q^{n-t}-1}\geq 1
\end{align*}}
from $\frac{q(q^{t-1}-1)}{q^{n-t}-1}\leq\frac{q(q^{t-1}-1)}{q^t-1}=1-\frac{q-1}{q^t-1}$. Thus $g_5(k_1,k_2,n,t)g_2(k_2,n,t)<g_3(k_1,k_2,n,t)$ by (\ref{v-ineq-11}).

Assume that $k_2=2t$. Since
{\small\begin{align*}
q^{t(k_1-t-2)}{n-2t-1\brack k_1-t-2}{n-t-2\brack k_1-t-2}^{-1}=\prod_{i=0}^{k_1-t-3}\frac{q^t(q^{n-2t-1-i}-1)}{q^{n-t-2-i}-1}\geq 1
\end{align*}}
and
{\small\begin{align*}
\frac{q^{n-t}-1}{q^2-1}\cdot\frac{q^{k_1-t-1}-1}{q^{n-t-1}-1}\cdot q^{-k_1+t+3}\geq  q^{n-t-2}\cdot q^{k_1-n-1} \cdot q^{-k_1+t+3} = 1,
\end{align*}}
by (\ref{v-ineq-11}), we have
{\small\begin{align*}
&\left(g_3(k_1,k_2,n,t)-g_5(k_1,k_2,n,t)g_2(k_1,k_2,n,t)\right){t+1\brack 1}^{-1}{n-t-1\brack k_1-t-1}^{-1}{n-t\brack t}^{-1}\\
>&\frac{q^{t}(q^{n-2t}-1)}{q^{n-t}-1}-1+q^{(t-1)(k_1-t-2)+1}\cdot\frac{q^{t}-1}{q^2-1}{n-2t-1\brack k_1-t-2}{n-t-1\brack k_1-t-1}^{-1}\\
=&-\frac{q^t-1}{q^{n-t}-1}+\frac{q^t-1}{q^2-1}\cdot\frac{q^{k_1-t-1}-1}{q^{n-t-1}-1}\cdot q^{-k_1+t+3}\cdot q^{t(k_1-t-2)}{n-2t-1\brack k_1-t-2}{n-t-2\brack k_1-t-2}^{-1}\\
\geq&\frac{q^t-1}{q^{n-t}-1}\left(-1 + \frac{q^{n-t}-1}{q^2-1}\cdot\frac{q^{k_1-t-1}-1}{q^{n-t-1}-1}\cdot q^{-k_1+t+3}\right)\geq 0,
\end{align*}}
implying that $g_5(k_1,k_2,n,t)g_2(k_2,n,t)<g_3(k_1,k_2,n,t)$ holds.

Therefore, the required result follows.  $\qed$
\begin{lemma}\label{v-inequ-3-1}
If $k_2=t+1,$ and $(k_1,k_2)\neq (2,2)$, $(3,2)$ or $(4,2)$, then $g_1(k_1,k_2,n,t)g_2(k_2,n,t)<g_3(k_1,k_2,n,t).$
\end{lemma}
\proof
By Lemma~\ref{lem1-1-1} (i), we have
{\small\begin{align*}
{n-t\brack k_1-t}-q^{2(k_1-t)}{n-t-2\brack k_1-t}={n-t-1\brack k_1-t-1}+q^{k_1-t}{n-t-2\brack k_1-t-1}={n-t-1\brack k_1-t-1}\left(1+\frac{q^{n-t}-q^{k_1-t}}{q^{n-t-1}-1}\right)
\end{align*}}
and ${n-t\brack 1}=q{n-t-1\brack 1}+1$, implying that
{\small\begin{align*}
g_1(k_1,k_2,n,t)g_2(k_2,n,t)&=\left({n-t\brack k_1-t}-q^{2(k_1-t)}{n-t-2\brack k_1-t}\right)\left({n-t\brack 1}+q^2{t\brack 1}\right)\\
&={n-t-1\brack k_1-t-1}\left(1+\frac{q^{n-t}-q^{k_1-t}}{q^{n-t-1}-1}\right)\left(q{n-t-1\brack 1}+1+q^2{t\brack 1}\right).
\end{align*}}
Observe that
{\small\begin{align*}
g_3(k_1,k_2,n,t)={n-t-1\brack k_1-t-1}\left(q{t+1\brack 1}{n-t-1\brack 1}+1\right).
\end{align*}}

Set
{\small\begin{align*}
g^\prime(k_1,k_2,n,t)=\left(g_3(k_1,k_2,n,t)-g_1(k_1,k_2,n,t)g_2(k_2,n,t)\right){n-t-1\brack k_1-t-1}^{-1}.
\end{align*}}
Then
{\small\begin{align*}
g^\prime(k_1,k_2,n,t)&=q{t+1\brack 1}{n-t-1\brack 1}+1-\left(1+\frac{q^{n-t}-q^{k_1-t}}{q^{n-t-1}-1}\right)\left(q{n-t-1\brack 1}+1+q^2{t\brack 1}\right)\\
&=q^4{t\brack 1}{n-t-3\brack 1}+q^3{t\brack 1}-q^{k_1-t+1}{n-k_1\brack 1}-\frac{q^{n-t}-q^{k_1-t}}{q^{n-t-1}-1}\cdot \left(1+q^2{t\brack 1}\right)
\end{align*}}
due to
{\small\begin{align*}
&q{t+1\brack 1}{n-t-1\brack 1}+1-\left(q{n-t-1\brack 1}+1+q^2{t\brack 1}\right)\\
=&q\left(q{t\brack 1}+1\right)\left(q{n-t-2\brack 1}+1\right)+1-\left(q^2{n-t-2\brack 1}+q+1+q^2{t\brack 1}\right)\\
=&q^3{t\brack 1}{n-t-2\brack 1}=q^4{t\brack 1}{n-t-3\brack 1}+q^3{t\brack 1}.
\end{align*}}
Observe that $k_1\geq 5$ if $t=1$. Since ${t\brack 1}\geq q^{t-1}$ and $\frac{q^{n-t}-q^{k_1-t}}{q^{n-t-1}-1}\leq q$, we have
{\small\begin{align*}
g^\prime(k_1,k_2,n,t)\geq &q^{4+(t-1)}{n-t-3\brack 1}-q^{k_1-t+1}{n-k_1\brack 1}-q\\
=&q(q-1)^{-1}(q^{n-1}-q^{n-t}-q^{t+2}+q^{k_1-t}-q+1)>0.
\end{align*}}
Hence, the desired result follows.  $\qed$
\begin{lemma}\label{v-inequ-6}
We have $$\max\{f_3(k_1,k_2,n,t)f_4(k_1,k_2,n,t),\ f_3(k_2,k_1,n,t)f_4(k_2,k_1,n,t)\}<g_4(k_1,k_2,n,t)g_2(k_2,n,t).$$	
\end{lemma}
\proof Assume that $(k,\ell)\in\{(k_1,k_2),\ (k_2,k_1)\}$. Since $n\geq k_1+k_2+t+3=k+\ell+t+3$, we have
{\small\begin{align*}
	q^{\ell-t+1}{t\brack 1}{k-t\brack 1}{n-t-2\brack \ell-t-1}{n-t\brack \ell-t}^{-1}=&q^{\ell-t+1}\cdot {t\brack 1}{k-t\brack 1}\cdot\frac{(q^{\ell-t}-1)(q^{n-\ell}-1)}{(q^{n-t}-1)(q^{n-t-1}-1)}\\
\leq& q^{2}\cdot {t\brack 1}{k-t\brack 1}\cdot\frac{q^{\ell-t}-1}{q^{n-t}-1}\cdot\frac{q^{n-t-1}-q^{\ell-t-1}}{q^{n-t-1}-1}\\	
\leq& q^{2+t+k-t+\ell-n}\leq q^{-1-t}\leq q^{-2},
\end{align*}}
implying that $f_4(k,\ell,n,t)\leq(1+q^{-2}){n-t\brack \ell-t}$. From $g_2(k_2,n,t)\geq {n-t\brack k_2-t}$, observe that
{\small\begin{align}
&g_4(k_1,k_2,n,t)g_2(k_2,n,t)-f_3(k,\ell,n,t)f_4(k,\ell,n,t)\nonumber\\	
	> &g_4(k_1,k_2,n,t){n-t\brack k_2-t}-f_3(k,\ell,n,t)\cdot(1+q^{-2}){n-t\brack \ell-t}.\label{ineq1}
\end{align}}
When $(k,\ell)=(k_1, k_2)$, since
{\small\begin{align*}
{k_2-t+1\brack 1}-(1+q^{-2}){k_2-t\brack 1}=q^{k_2-t}-q^{-2}{k_2-t\brack 1}\geq q^{k_2-t}-q^{-2+k_2-t}>q^{k_2-t-1},
\end{align*}}
by (\ref{ineq1}), we have
{\small\begin{align*}
	&\left(g_4(k_1,k_2,n,t)g_2(k_2,n,t)-f_3(k_1,k_2,n,t)f_4(k_1,k_2,n,t)\right){n-t-2\brack k_1-t-2}^{-1}{n-t\brack k_2-t}^{-1}\\	
	> &{k_2-t+1\brack 1}\frac{q^{n-t-1}-1}{q^{k_1-t-1}-1}-q{k_2-t+1\brack 2}-\left({k_2-t\brack 1}\frac{q^{n-t-1}-1}{q^{k_1-t-1}-1}+q^{2(k_2-t)}\right)\cdot(1+q^{-2})\\
=&\left({k_2-t+1\brack 1}-(1+q^{-2}){k_2-t\brack 1}\right)\frac{q^{n-t-1}-1}{q^{k_1-t-1}-1}-q{k_2-t+1\brack 2}-(1+q^{-2})q^{2(k_2-t)}\\
>&q^{k_2-t-1}\cdot q^{n-k_1}-q^{2(k_2-t)+1}-(1+q^{-2})q^{2(k_2-t)}\\
>&q^{2k_2+2}-2q^{2(k_2-t)+1}\geq 0,
\end{align*}}
implying that $f_3(k_1,k_2,n,t)f_4(k_1,k_2,n,t)<g_4(k_1,k_2,n,t)g_2(k_2,n,t)$ holds. When $(k,\ell)=(k_2, k_1)$, by (\ref{ineq1}), we have
{\small \begin{align*}
	&\left(g_4(k_1,k_2,n,t)g_2(k_2,n,t)-f_3(k_2,k_1,n,t)f_4(k_2,k_1,n,t)\right){n-t-1\brack k_1-t-1}^{-1}{n-t-1\brack k_2-t-1}^{-1}{n-t\brack 1}^{-1}\\		 >&\frac{q^{k_2-t+1}-1}{q^{k_2-t}-1}-\frac{q(q^{k_2-t+1}-1)(q^{k_1-t-1}-1)}{(q^{2}-1)(q^{n-t-1}-1)}-(1+q^{-2})\left(1+\frac{q^{2(k_1-t)}(q-1)(q^{k_2-t-1}-1)}{(q^{k_1-t}-1)(q^{n-t-1}-1)}\right)\\
	 >&q-(1+q^{-2})-q^{1+(k_2-t)+(k_1-n)}-(1+q^{-2})q^{2(k_1-t)+(1-k_1+t)+(k_2-n)}\\
	\geq&q-1-q^{-2}-(2+q^{-2})\cdot q^{-2t-2}>0,
\end{align*}}
and $f_3(k_2,k_1,n,t)f_4(k_2,k_1,n,t)<g_4(k_1,k_2,n,t)g_2(k_2,n,t)$ holds.
$\qed$

\begin{lemma}\label{v-inequ-5}
	If $k_1>k_2,$ then $g_1(k_2,k_1,n,t)g_2(k_1,n,t)<g_1(k_1,k_2,n,t)g_2(k_2,n,t)$.
\end{lemma}
\proof From Lemma~\ref{lem-upp-1}, we have $g_4(k_1,k_2,n,t)\leq g_1(k_1,k_2,n,t)$ and $g_1(k_2,k_1,n,t)\leq g_5(k_2,k_1,n,t)\leq{k_1-t+1\brack 1}{n-t-1\brack k_2-t-1}.$ Thus
{\small\begin{align*}
	&g_1(k_1,k_2,n,t)g_2(k_2,n,t)-g_1(k_2,k_1,n,t)g_2(k_1,n,t)\\
	=&\left(g_1(k_1,k_2,n,t){n-t\brack k_2-t}-g_1(k_2,k_1,n,t){n-t\brack k_1-t}\right)\\
	&+\left(g_1(k_1,k_2,n,t)\cdot q^{k_2+1-t}{t\brack 1}-g_1(k_2,k_1,n,t)\cdot q^{k_1+1-t}{t\brack 1}\right)\\
	\geq&\left(q^{(k_1+1-t)(k_2-t)}{n-k_1-1\brack k_2-t}{n-t\brack k_1-t}-q^{(k_2+1-t)(k_1-t)}{n-k_2-1\brack k_1-t}{n-t\brack k_2-t}\right)\\
	&+q^{k_2+1-t}{t\brack 1}\left(g_4(k_1,k_2,n,t) -{k_1-t+1\brack 1}{n-t-1\brack k_2-t-1}\cdot q^{k_1-k_2}\right).
\end{align*}}
Since $n\geq k_1+k_2+t+3$, we have
{\small\begin{align*}
	\frac{q^{(k_1+1-t)(k_2-t)}{n-k_1-1\brack k_2-t}{n-t\brack k_1-t}}{q^{(k_2+1-t)(k_1-t)}{n-k_2-1\brack k_1-t}{n-t\brack k_2-t}}=&\frac{\prod_{i=0}^{k_2-t-1}(q^{n-k_1-1-i}-1)\cdot \prod_{j=0}^{k_1-t-1}(q^{n-t-j}-1)}{q^{k_1-k_2}\prod_{i=0}^{k_1-t-1}(q^{n-k_2-1-i}-1)\cdot \prod_{j=0}^{k_2-t-1}(q^{n-t-j}-1)}\\
	=&\frac{ \prod_{j=k_2-t}^{k_1-t-1}(q^{n-t-j}-1)}{\prod_{i=0}^{k_1-k_2-1}(q^{n-k_2-i}-q)}>1
\end{align*}}
and
{\small\begin{align*}
	&\left(g_4(k_1,k_2,n,t)-{k_1-t+1\brack 1}{n-t-1\brack k_2-t-1}\cdot q^{k_1-k_2}\right){k_2-t+1\brack 1}^{-1}{n-t-2\brack k_1-t-2}^{-1}\\
=&\frac{q^{n-t-1}-1}{q^{k_1-t-1}-1}-\frac{q(q^{k_2-t}-1)}{q^2-1}-\frac{q^{k_1-k_2}(q^{k_1-t+1}-1)(q^{n-t-1}-1)}{(q^{k_2-t+1}-1)(q^{n-k_2}-1)}\cdot \prod_{j=0}^{k_1-k_2-2}\frac{q^{k_1-t-2-j}-1}{q^{n-k_2-1-j}-1}\\
\geq&q^{n-k_1}-q^{k_2-t}-q^{(k_1-k_2)+(k_1-k_2+1)+(k_2-t)+(k_1-k_2-1)(k_1+k_2-t-1-n)}\\
\geq&q^{k_2+t+3}-q^{k_2-t}-q^{2k_1-k_2+1-t+(k_1-k_2-1)(-2t-4)}\\
=&q^{k_2+t+3}-q^{k_2-t}-q^{k_2-t+3+(k_1-k_2-1)(-2t-2)}>0.
\end{align*}}
Hence the required result follows. $\qed$

\begin{lemma}\label{v-inequ-1}
If $k_1>k_2$, then $g_3(k_2,k_1,n,t)<g_3(k_1,k_2,n,t).$	
\end{lemma}
\proof Observe that
{\small\begin{align*}
	\frac{g_3(k_1,k_2,n,t)-g_3(k_2,k_1,n,t)}{{n-t-1\brack k_1-t-1}{n-t-1\brack k_2-t-1}{t+1\brack 1}}=& \frac{q^{n-t}-q^{k_2-t}}{q^{k_2-t}-1}-\frac{q^{n-t}-q^{k_1-t}}{q^{k_1-t}-1}\\
	=&\frac{(q^{n-t}-1)(q^{k_1-t}-q^{k_2-t})}{(q^{k_2-t}-1)(q^{k_1-t}-1)}>0,
\end{align*}}
and the required result follows.   $\qed$

\begin{lemma}\label{v-inequ-8}
We have ${n-t\brack k_2-t}f_2(t+2,k_2,k_1,n,t)<g_4(k_1,k_2,n,t)g_2(k_2,n,t).$	
\end{lemma}
\proof It suffices to prove that $f_2(t+2, k_2, k_1,n,t)<g_4(k_1,k_2,n,t)$.
Since $\frac{1}{q+1}+{t+2\brack 2}<1+{t+2\brack 2}\leq q^{2(t+1)}$, we have
{\small\begin{align*}
	&\left(g_4(k_1,k_2,n,t)-f_2(t+2,k_2,k_1,n,t)\right)\cdot{k_2-t+1\brack 1}^{-1}{n-t-2\brack k_1-t-2}^{-1}\\
	=&\frac{q^{n-t-1}-1}{q^{k_1-t-1}-1}-\frac{q(q^{k_2-t}-1)}{(q+1)(q-1)}-{t+2\brack 2}{k_2-t+1\brack 1}\\
>&\frac{q^{n-t-1}-1}{q^{k_1-t-1}-1}-\frac{q^{k_2-t+1}}{q+1}-q^{k_2-t+1}{t+2\brack 2}\\
>&q^{n-k_1}-q^{2(t+1)+k_2-t+1}\geq 0
\end{align*}}
due to $n\geq k_1+k_2+t+3$. Therefore the desired result follows. $\qed$

\begin{lemma}\label{v-inequ-9}
The following hold.
\begin{enumerate}[{\rm(i)}]
\item If $k_2\leq 2t$, then	
${n-t\brack k_1-t}f_2(t+2,k_1,k_2,n,t)<g_3(k_1,k_2,n,t).$
\item If $k_2\geq 2t+1$, then	
${n-t\brack k_1-t}f_2(t+2,k_1,k_2,n,t)<g_4(k_1,k_2,n,t)g_2(k_2,n,t).$
\end{enumerate}
\end{lemma}
\proof (i) Since
{\small\begin{align*}
	&(q^{n-t}-q^{k_2-t})(q^2-1)-(q^{n-t}-1)(q^{t+2}-1)q^{-t-1}\\
	=&q^{n-t+2}-q^{n-t}-q^{k_2-t+2}+q^{k_2-t}-q^{n-t+1}+q^{n-2t-1}+q-q^{-t-1}\\
	=&q^{n-t+2}-q^{n-t+1}-q^{n-t}+q^{k_2-t}+q^{n-2t-1}-q^{k_2-t+2}+q-q^{-t-1}>0,
\end{align*}}
we have $\frac{q^{k_2-t}(q^{n-k_2}-1)(q^2-1)}{(q^{n-t}-1)(q^{t+2}-1)}>q^{-t-1}$. Hence
{\small\begin{align*}
\frac{g_3(k_1,k_2,n,t)}{{n-t\brack k_1-t}f_2(t+2,k_1,k_2,n,t)}>&\frac{{n-t-1\brack k_1-t-1}\cdot q^{k_2-t}{t+1\brack 1}{n-t-1\brack k_2-t}}{{n-t\brack k_1-t}{t+2\brack 2}{k_1-t+1\brack 1}^2{n-t-2\brack k_2-t-2}}\\
=&\frac{(q^{k_1-t}-1)q^{k_2-t}(q^2-1)(q^{n-t-1}-1)(q^{n-k_2}-1)}{(q^{n-t}-1)(q^{t+2}-1)(q^{k_2-t}-1)(q^{k_2-t-1}-1)}\cdot{k_1-t+1\brack 1}^{-2}\\
>&\frac{q^{k_1-t}-1}{q^{k_2-t}-1}\cdot\frac{q^{n-t-1}-1}{q^{k_2-t-1}-1}\cdot q^{-t-1}\cdot {k_1-t+1\brack 1}^{-2}\\
>&q^{(k_1-k_2)+(n-k_2)+(-t-1)-2(k_1-t+1)}=q^{n-2k_2-k_1+t-3}\geq 1
\end{align*}}
due to $n\geq k_1+k_2+t+3$ and $k_2\leq 2t,$ and (i) holds.

\medskip

(ii) By $k_1\geq k_2\geq 2t+1$ and Lemma~\ref{lem1-1-1} (ii) and (iii), we have
{\small\begin{align*}
	\frac{{n-t\brack k_1-t}{t+2\brack 2}{k_1-t+1\brack 1}^2{n-t-2\brack k_2-t-2}}{{k_2-t+1\brack 2}{n-t-2\brack k_1-t-2}{n-t\brack k_2-t}}=&\frac{q^{t+2}-1}{q^{k_2-t+1}-1}\cdot\frac{q^{t+1}-1}{q^{k_1-t}-1}\cdot\frac{q^{k_2-t-1}-1}{q^{k_1-t-1}-1}\cdot{k_1-t+1\brack 1}^2\\
	<&q^{(2t-k_2+1)+(2t-k_1+1)+(k_2-k_1)+2(k_1-t+1)}=q^{2t+4}
\end{align*}}
implying that
{\small\begin{align*}
{n-t\brack k_1-t}f_2(t+2,k_1,k_2,n,t)=&{n-t\brack k_1-t}{t+2\brack 2}{k_1-t+1\brack 1}^2{n-t-2\brack k_2-t-2}\\
<&q^{2t+4}{k_2-t+1 \brack 2}{n-t-2 \brack k_1-t-2}{n-t \brack k_2-t}.
\end{align*}}
Then
{\small\begin{align*}
	\frac{g_4(k_1,k_2,n,t){n-t\brack k_2-t}-{n-t\brack k_1-t}f_2(t+2,k_1,k_2,n,t)}{{k_2-t+1 \brack 2}{n-t-2 \brack k_1-t-2}{n-t \brack k_2-t}}
	>&\frac{(q+1)(q-1)(q^{n-t-1}-1)}{(q^{k_2-t}-1)(q^{k_1-t-1}-1)}-q-q^{2t+4}\\
>&(q+1)q^{(-k_2+t)+(n-k_1)}-q-q^{2t+4}>0
\end{align*}}
from $n\geq k_1+k_2+t+3,$ and (ii) holds.   $\qed$
\begin{lemma}\label{v-inequ-4}
	We have $f_1(k_1,k_2,n,t)f_1(k_2,k_1,n,t)<g_4(k_1,k_2,n,t)g_2(k_2,n,t)$.
\end{lemma}
\proof Since $g_2(k_2,n,t)>{n-t\brack k_2-t}$, by Lemma~\ref{lem1-1-1} (ii), we have
{\small\begin{align*}
&\frac{g_4(k_1,k_2,n,t)g_2(k_2,n,t)-f_1(k_1,k_2,n,t)f_1(k_2,k_1,n,t)}{{k_1-t+1\brack 1}{k_2-t+1\brack 1}{n-t-1\brack k_1-t-1}{n-t-1\brack k_2-t-1}}\\	
>&\frac{(q^{n-t}-1)(q-1)}{(q^{k_2-t}-1)(q^{k_1-t+1}-1)}- \frac{q(q^{n-t}-1)(q^{k_1-t-1}-1)}{(q+1)(q^{n-t-1}-1)(q^{k_1-t+1}-1)}-{t+1\brack 1}^2\\
\geq&q^{n-k_2-k_1+t-1}-1-{t+1\brack 1}^2\\
=&(q-1)^{-2}((q^{n-k_2-k_1+t-1}-1)(q^2-2q)+q^{n-k_2-k_1+t-1}-1-q^{2t+2}+2q^{t+1}-1)>0
\end{align*}}
due to $n\geq k_1+k_2+t+3$ and $t\geq 1.$ Therefore, the required result follows. $\qed$

\section*{Acknowledgement}

M. Lu is supported by the National Natural Science Foundation of China (12171272), B. Lv is supported by National Natural Science Foundation of China (12071039, 12131011), K. Wang is supported by the National Key R\&D Program of China (No. 2020YFA0712900) and National Natural Science Foundation of China (12071039, 12131011).


\end{document}